\documentclass[letterpaper, 11pt]{article}

\usepackage[english]{babel}
\usepackage[left=1 in, right=1 in, top=1.5 in, bottom=1.5 in]{geometry}
\usepackage{amsmath}
\usepackage{amssymb}
\usepackage{mathtools}
\usepackage{mathrsfs}
\usepackage[utf8]{inputenc}
\usepackage[dvipsnames]{xcolor}
\usepackage{amsthm}
\usepackage{tikz, calc}
\usetikzlibrary{graphs,arrows.meta,shapes.geometric}
\usepackage[shortlabels]{enumitem}
\usepackage{float}
\usepackage{graphicx}
\usepackage{caption}
\usepackage{subcaption}
\usepackage{tabularx}
\usepackage{booktabs}
\usepackage{hyperref}
\usepackage{cleveref}
 
\newcommand{\cupdot}{\mathbin{\mathaccent\cdot\cup}}

\makeatletter
\DeclareRobustCommand{\rvdots}{%
  \vbox{
    \baselineskip4\p@\lineskiplimit\z@
    \kern-\p@
    \hbox{.}\hbox{.}\hbox{.}
  }}
\makeatother

\theoremstyle{definition}
\newtheorem{theorem}{Theorem}[section]
\newtheorem{definition}[theorem]{Definition}
\newtheorem{proposition}[theorem]{Proposition}
\newtheorem{corollary}[theorem]{Corollary}
\newtheorem{lemma}[theorem]{Lemma}
\newtheorem{example}[theorem]{Example}
\newtheorem{remark}[theorem]{Remark}

\newcommand{\bsy}{\boldsymbol}

\newcommand{\cpf}{\mathcal{P}\mkern-3mu\mathcal{F}}
\newcommand{\cMPF}{\mathcal{MP}\mkern-3mu\mathcal{F}} 

\newcommand{\indeg}{\mathrm{indeg}}
\newcommand{\outdeg}{\mathrm{outdeg}}

\title{($\mathfrak{S}_p \times \mathfrak{S}_q$)-Invariant Graphical Parking Functions}

\author{Lauren Snider\footnote{laurenleighsnider@gmail.com} and Catherine Yan\footnote{huafei-yan@tamu.edu, C.H.~Yan was supported  by the Simons Collaboration Grant for Mathematics 704276. } 
\\
Department of Mathematics, Texas A\&M University, \\
College Station, TX 77843}

\date{}

\begin{document}

\maketitle

\begin{abstract}
Graphical parking functions, or $G$-parking functions, are a generalization of classical parking functions that depend on a connected multigraph $G$ with a distinguished root vertex. Gaydarov and Hopkins established a connection  between $G$-parking functions and a vector-dependent generalization of parking functions known as $\bsy{u}$-parking functions. The central component of their result was the classification of all graphs $G$ for which the set of  $G$-parking functions is invariant under the action of the symmetric group $\mathfrak{S}_n$, where $n+1$ is the order of $G$. In this work, we present a  higher dimensional analogue of Gaydarov and Hopkins' results by characterizing the intersection between $G$-parking functions and 2-dimensional $\bsy{U}$-parking functions, which are pairs of integer sequences whose order statistics are bounded by certain weights along lattice paths in the plane. Our key result is a complete charaterization of all $G$ for which the set of $G$-parking functions is invariant under the action of $\mathfrak{S}_p \times \mathfrak{S}_q$, where $p+q+1$ is the order of $G$.  
\end{abstract}



\section{Introduction}

Classical parking functions were first considered by Konheim and Weiss in 1966 as an illustrative device in their work on hashing functions, with their original definition in terms of a parking problem \cite{KW66}. Consider a one-way street with $n$ labeled parking spaces and the same number of drivers, each having a preferred parking space. Drivers enter the street one-by-one and attempt to park in their preferred space or, if unavailable, the next available space down the street (if one exists). The sequence of driver preferences is said to be a \textit{(classical) parking function} if all drivers can park according to this process. 
Let $\mathbb{N}=\{0, 1, 2, \dots\}$ be the set of non-negative integers. 
The following is an equivalent, albeit more formal, definition of parking function. 

\begin{definition}
\label{def:PF2}
A \textit{parking function of length n} is a sequence $\bsy{a} = (a_1,a_2, \ldots, a_n) \in \mathbb{N}^n$ whose non-decreasing rearrangement $a_{(1)} \leq a_{(2)} \leq \cdots \leq a_{(n)}$ satisfies $a_{(i)} < i$ for each $i$. The non-decreasing sequence $(a_{(1)}, a_{(2)}, \ldots, a_{(n)})$ is the \textit{order statistics of} $\bsy{a}$, with $a_{(i)}$ the \textit{i-th order statistic of} $\bsy{a}$.
\end{definition}

It is evident from the above definition that any permutation of the entries of a parking function must also be a parking function, that is, the symmetric group $\mathfrak{S}_n$ on $n$ elements acts on the set $\cpf_n$ of length-$n$ parking functions via permutations of coordinates.

Since their introduction, parking functions have become a fertile source of combinatorial research, leading to numerous generalizations and surprising connections to other concepts and disciplines, including hyperplane arrangements \cite{S98}, interpolation theory \cite{KY03}, combinatorial theory of McDonnald polynomials \cite{Haglund08},   the Abelian sandpile model \cite{CP02,PS04}, and many others; see the extensive survey \cite{Y15} for more on the combinatorial theory of parking functions. 

One of the most basic generalizations of parking functions is the \textit{vector parking function}, or \textit{$\bsy{u}$-parking function} if the vector $\bsy{u}$ is specified.

\begin{definition}
\label{def:vector_PF}
Let $\bsy{u} = (u_1, u_2, \ldots, u_n)$ be a non-decreasing sequence of positive integers. A sequence $\bsy{a} = (a_1, a_2, \ldots, a_n) \in \mathbb{N}^n$ is a \textit{$\bsy{u}$-parking function} of length $n$ if the order statistics of $\bsy{a}$ satisfy $a_{(i)} < u_i$ for each $i$.
\end{definition}

The set of $\bsy{u}$-parking functions of length $n$ will be denoted by $\cpf_n(\bsy{u})$. Notably, when $\bsy{u}$ is the arithmetic sequence $(1,2,\ldots,n)$, the notions of $\bsy{u}$-parking functions and classical parking functions are consistent. 

Another generalization of parking functions is given by \textit{graphical parking functions}, or \textit{$G$-parking functions} for a designated graph $G$. Though their name and first formal definition in the literature are attributed to Postnikov and Shapiro in \cite{PS04}, their idea can be traced back further in the context of the Abelian sandpile model (see \cite{BTW87} or \cite{D90}). Postnikov and Shapiro noted that $G$-parking functions are in essence dual to the \textit{recurrent configurations} of the Abelian sandpile model \cite[Lemma 13.6]{PS04}.

As suggested by their name, the set of $G$-parking functions, denoted by $\cpf(G)$, depends on a choice of graph $G$. Gaydarov and Hopkins clarified the relationship between graphical parking functions and vector parking functions in \cite{GH16} by characterizing the structure of all graphs $G$ and all vectors $\bsy{u}$ for which $\cpf(G) = \cpf(\bsy{u})$. The cornerstone of their proof was the classification of all graphs $G$ whose set of $G$-parking functions is invariant under the action of $\mathfrak{S}_n$, where $n+1$ is the order of $G$. Notably, they found exactly three graph structures having the specified invariance property: trees, cycles, and complete graphs (each having certain requirements for edge-weights).

The intent of the present paper is to extend Gaydarov and Hopkins' results to two dimensions. Doing so requires the introduction of a third generalization of parking functions: the multidimensional $\bsy{U}$-parking functions. 
Proposed by Khare, Lorentz, and Yan in their study of  multivariate Gon\v{c}arov polynomials \cite{KLY14}, $k$-dimensional $\bsy{U}$-parking functions are $k$-tuples of integer sequences whose order statistics are bounded above by edge-weights on lattice paths, with these edge-weights given by a $k$-dimensional array $\bsy{U}$. Each element of $\bsy{U}$ is a node in $\mathbb{N}^k$, and 
we refer to $\bsy{U}$ as the \emph{weight set} of the $k$-dimensional parking functions. 
When $k=1$, meaning the set of nodes in  $\bsy{U}$ 
form a vector, the notion of a $\bsy{U}$-parking function is in accordance with the definition of a vector parking function. 

In this paper, we will restrict our attention to the set $\cpf_{p,q}^{(2)}(\bsy{U})$ of 2-dimensional $\bsy{U}$-parking functions in order to address the question: for what graphs $G$ and weight sets $\bsy{U}$ do we have $\cpf(G) = \cpf_{p,q}^{(2)}(\bsy{U})$? Similar to the method of Gaydarov and Hopkins' proof, the key ingredient is to   classify all graphs $G$ whose set of $G$-parking functions is invariant under the action of the product of symmetric groups $\mathfrak{S}_p \times \mathfrak{S}_q$, where $p+q+1$ is the order of $G$.

The remainder of the paper is organized as follows. In Section 2, we review the definition of $G$-parking functions and Gaydarov and Hopkins' classification, along with an important result concerning acyclic orientations of $G$. Section 3 provides a brief overview of 2-dimensional $\bsy{U}$-parking functions. In Section 4 
We focus on the special case of affine  $\bsy{U}$, that is, the nodes of $\bsy{U}$ can be written as $(u_{i,j},v_{i,j})^T = A(i,j)^T + (c_1,c_2)^T$, with $A$ a $2 \times 2$ integer matrix and $c_1,c_2 \in \mathbb{N}$.  We determine the affine weight set $\bsy{U}$ for which there exist graphs $G$ satisfying $\cpf(G)=\cpf_{p,q}^{(2)}(\bsy{U})$. 
In  \Cref{sec:classification}, we classify all graphs whose corresponding set of graphical parking functions is invariant under the action of $\mathfrak{S}_p \times \mathfrak{S}_q$, where $p+q+1$ is the order of $G$. Finally, in \Cref{sec:finish}, for each graph $G$ whose $G$-parking functions are $(\mathfrak{S}_p \times \mathfrak{S}_q)$-invariant, we identify all the 
weight sets $\bsy{U}$  for which $\cpf_{p,q}^{(2)}(\bsy{U})$ and $\cpf(G)$ coincide, thereby completing our characterization of the overlap between graphical parking functions and 2-dimensional $\bsy{U}$-parking functions.


\section{Graphical parking functions}
\label{section 2}

The notion of $G$-parking function was first introduced by 
Postnikov and Shapiro in \cite{PS04}. In this section, we review its definition for un-directed graphs, along with two relevant results: the correspondence between maximal $G$-parking functions and certain acyclic orientations of $G$, and the relationship between $G$-parking functions and vector parking functions.

For the rest of this paper, any graph $G$ is assumed to be an undirected, connected, loopless multigraph having a distinguished root vertex unless otherwise specified. For our purposes, we can equivalently view $G$ as a simple graph (i.e., having no loops or parallel edges) with vertex set $V(G)$, edge set $E(G)$, and an associated edge-weight function $wt_G: E(G) \to \mathbb{Z}^+$, where $wt_G(\{i,j\})$ equals the number of parallel edges between vertices $i$ and $j$ in the multigraph $G$. If vertices $i$ and $j$ are not adjacent in $G$, we may sometimes write $wt_G(\{i,j\})=0$ to indicate $\{i,j\} \notin E(G)$. We will always let $V(G) = [n]_0 = \{0,1,2,\ldots,n\}$ for some $n \in \mathbb{Z}^+$, with $0$ the unique root vertex.

If $U$ is a subset of $V(G)$ and $i \in U$, we write $d_U(i)$ to denote the sum $\sum_{j \not \in U} wt_G(\{i,j\})$.
Let $[n]=\{1, 2, \dots, n\}$. 
We are now able to define the notion of a graphical parking function, or $G$-parking function.

\begin{definition}
\label{def:G-PF}
Let $G$ be a graph with $V(G) = [n]_0$. A $G$-\textit{parking function} is a function $f: [n]\to \mathbb{N}$ 
such that, for each non-empty subset $U \subseteq [n]$, there exists a vertex $i \in U$ having $f(i) < d_U(i)$. A $G$-parking function is written as the sequence $(f(1), f(2), \dots, f(n))$. 
\end{definition}

We will denote the set of $G$-parking functions by $\cpf(G)$. 
Clearly $\cpf(G)\neq \emptyset$ since $f(i)=0$ for all $i \in [n]$ is always a $G$-parking function. 
When $G$ is the complete graph $K_{n+1}$ having vertex set $[n]_0$ and $wt_G(e) = 1$ for all $e \in E(G)$, the set of $G$-parking functions is precisely the set of classical parking functions. 

There is a natural partial order $\preceq$ on $\mathbb{N}^n$ given by: $(a_1,a_2,\ldots,a_n) \preceq (b_1,b_2,\ldots,b_n)$ if and only if $a_i \leq b_i$ for all $i \in [n]$. 
Denote  by $\cMPF(G)$ the set of maximal $G$-parking functions under the order $\preceq$. It is easily seen that $\bsy{a} \in \cpf(G)$ if and only if $\bsy{a} \preceq \bsy{b}$ for some $\bsy{b} \in \cMPF(G)$. Hence, finding all $G$-parking functions for a given graph $G$ can be reduced to the problem of computing the set $\cMPF(G)$.

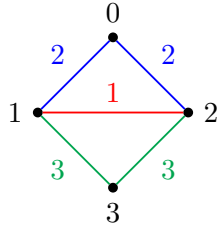
\begin{figure}[H]
\begin{center}
\begin{tikzpicture}[line width=.7pt,scale=1]
\tikzstyle{vertex}=[draw,circle,fill=black,minimum size=3pt,inner sep=0pt]

\path (0,0) node[vertex, label=above:0] (v0) {};
\path (-1,-1) node[vertex, label=left:1] (v1) {};
\path (1, -1) node[vertex, label=right:2] (v2) {};
\path (0,-2) node[vertex, label=below:3] (v3) {};

\draw[blue] (v0) -- node [above left] {\textcolor{blue}{2}} (v1) ;
\draw[blue] (v0) -- node [above right] {\textcolor{blue}{2}} (v2);
\draw[red] (v1) -- node [above] {\textcolor{red}{1}} (v2);
\draw[Green] (v1) -- node [below left] {\textcolor{Green}{3}} (v3);
\draw[Green] (v2) -- node [below right] {\textcolor{Green}{3}} (v3);

\end{tikzpicture}
\caption{
In the graph $G$ above,
$\cMPF(G)=\{(5,1,2),(1,5,2), (2,1,5), (1,2,5)\}.$  
Note that all edges of $G$ having the same weight are represented by the same color.
}
\label{fig:G-PG_example}
\end{center}
\end{figure}

\subsection{Maximal $G$-parking functions and acyclic orientations}
\label{subsection 2.1}

Maximal $G$-parking functions can be represented  in terms of certain acyclic orientations of $G$. Here, an \textit{orientation} is a subset $\mathcal{O} \subseteq V(G) \times V(G)$ such that:
\begin{enumerate}
    \item $\{i,j\} \in E(G)$ if and only if $(i,j)\in\mathcal{O}$ or $(j,i)\in\mathcal{O}$, and
    \item $(i,j)\in\mathcal{O}$ implies $(j,i)\not\in\mathcal{O}$ for all $i,j \in V(G)$.
\end{enumerate}

The ordered pair $(i,j) \in \mathcal{O}$ means the edge $\{i,j\}$ is directed from $i$ to $j$. For an orientation $\mathcal{O}$ of $G$ and $i \in V(G)$, the \textit{indegree of i with respect to} $\mathcal{O}$ is given by $\indeg_{_\mathcal{O}}(i) = \sum_{(j,i)\in\mathcal{O}} w_G(\{i,j\})$, and the \textit{outdegree of i with respect to} $\mathcal{O}$ is given by $\mathrm{outdeg}_{_\mathcal{O}}(i) = \sum_{(i,j)\in\mathcal{O}} w_G(\{i,j\})$. A vertex $i \in V(G)$ is said to be a \textit{source of} $\mathcal{O}$ if $\indeg_{_\mathcal{O}}(i) = 0$, whereas $i$ is said to be a \textit{sink of} $\mathcal{O}$ if $\outdeg_{_\mathcal{O}}(i) = 0$. Also, $\mathcal{O}$ is \textit{acyclic} if there is no sequence $(v_1,v_2),(v_2,v_3),\ldots,(v_{m-1},v_m),(v_m,v_1)$ of directed edges in $\mathcal{O}$.

Adopting the notation of Gaydarov and Hopkins in \cite{GH16}, we will let $\mathcal{A}(G)$ denote the set of all acyclic orientations of $G$ having 0 as the unique source. 
The following theorem, relating  maximal graphical parking functions to the set $\mathcal{A}(G)$, is well-known and can be tracked back to a  result of Greene and Zaslavsky  \cite{GZ83}. See, for example, \cite[Theorem 3.1]{BCT10} and the references there.

\begin{proposition} 
\label{burningalgorithm}
Let $G$ be a multigraph with $V(G)=[n]_0$. There is a bijection between $\mathcal{A}(G)$ and $\cMPF(G)$, given by $\mathcal{O} \mapsto (\indeg_{_\mathcal{O}}(1)-1, \indeg_{_\mathcal{O}}(2)-1, \ldots, \indeg_{_\mathcal{O}}(n)-1)$.
\end{proposition}

We should remark that although some orientations of the multigraph $G$ are not realizable when $G$ is viewed as a simple graph with edge-weight function $wt_G$, this discrepancy turns out to be inconsequential since we will only be concerned with acyclic orientations in the set $\mathcal{A}(G)$, where both perspectives agree with each other.

\subsection{Overlap between $G$-parking functions and vector parking functions}
\label{subsection 2.2}

In \cite{GH16}, Gaydarov and Hopkins classified all vectors $\bsy{u}$ and graphs $G$ for which $\cpf(\bsy{u}) = \cpf(G)$. 
To state their result, we need some new graph terminology.  

Given a graph $G$ and its associated edge-weight function $wt_G: E(G) \to \mathbb{Z}^+$, we say:
\begin{itemize}
    \item $G$ is an \textit{$a$-tree} if $G$ is a tree with $wt_G(\{i,j\}) = a$ for all $\{i,j\} \in E(G)$. 
    \item $G$ is an \textit{$a$-cycle} if $G$ is a cycle with $wt_G(\{i,j\}) = a$ for all $\{i,j\} \in E(G)$. 
    \item $G$ is equal to $K_{n+1}^{a,b}$ if $G$ is the complete graph $K_{n+1}$ with vertex set $[n]_0$,  $wt_G(\{0,i\})=a$ for all  $i \in [n]$,  and $wt_G(\{i,j\})=b$ for all $i,j \in [n]$.  
\end{itemize}

\begin{theorem} \cite[Theorem 2.5]{GH16}
\label{thm:gaydarov and hopkins}
If $G$ is a graph such that $\cpf(G)$ is invariant under the action of $\mathfrak{S}_n$, then one of the following cases holds:
\begin{enumerate}
 \item[(i)] $\cpf((a,a,\ldots,a)) = \cpf(G)$, where $a \geq 1$ and $G$ is an $a$-tree;
 \item[(ii)] $\cpf((a,a,\ldots,a,2a)) = \cpf(G)$, where $a \geq 1$ and $G$ is an $a$-cycle;
 \item[(iii)] $\cpf((a,a+b,a+2b,\ldots,a+(n-1)b)) = \cpf(G)$, where $a,b,n \geq 1$ and $G$ is equal to $K_{n+1}^{a,b}$.
\end{enumerate}
Otherwise, if $\cpf(G)$ is not invariant under the action of $\mathfrak{S}_n$, then there is no $\bsy{u} \in (\mathbb{Z}^+)^n$ such that $\cpf(G) = \cpf(\bsy{u})$.
\end{theorem}

The essential ingredient of Gaydarov and Hopkins' proof 
is the classification of all graphs $G$ for which $\cMPF(G)$ is invariant under the action of $\mathfrak{S}_n$.
By working strictly with the set of maximal $G$-parking functions, they are able to exploit the bijective correspondence between 
$\cMPF(G)$ and $\mathcal{A}(G)$ of Proposition~\ref{burningalgorithm}. We will follow an analogous approach in our attempt to extend Theorem~\ref{thm:gaydarov and hopkins} to two dimensions.


\section{Two-dimensional $\bsy{U}$-parking functions}
\label{section 3}

In this section, we briefly review all definitions and associated notation relevant to the 2-dimensional $\bsy{U}$-parking function. These higher dimensional generalizations of parking functions were first introduced as combinatorial objects enumerated by the multivariate Gon\v{c}arov polynomials, a special class of interpolating polynomials considered by Khare, Lorentz, and Yan in \cite{KLY14}
and further explored in \cite{ASY21,LTY17,LY16,SY22}.  %
Higher dimensional $\bsy{U}$-parking functions and multivariate Gon\v{c}arov polynomials have applications in probability theory, epidemic models, and insurance risk processes \cite{LePh16,LePh23}.  In this paper, we restrict ourselves to the 2-dimensional case.

In the following discussion, assume $p,q \in \mathbb{N}$ are fixed.
Let $D_{p,q}$ be the directed graph whose vertices are the lattice points $\{(i,j) : 0 \leq i \leq p, 0 \leq j \leq q\}$ and whose edges consist of all north steps $N=(0,1)$ and east steps $E = (1,0)$ connecting its vertices. 
Let $ \bsy{U} \subset \mathbb{N}^2$ be a set of \textit{weights}  given by $\bsy{U} = \{(u_{i,j},v_{i,j}) : 0 \leq i \leq p, 0 \leq j \leq q \}$. We will denote by $D_{p,q}(\bsy{U})$ the weighted directed graph on $D_{p,q}$, where  an edge $e$ of  has an associated weight $w_U(e)$ given by:
\[ w_U(e) = \begin{cases} 
      u_{i,j} & \textnormal{if } e \textnormal{ is an \textit{E}-step from } (i,j) \textnormal{ to } (i+1,j), \\
      v_{i,j} & \textnormal{if } e \textnormal{ is an \textit{N}-step from } (i,j) \textnormal{ to } (i,j+1).
   \end{cases}
\]

We can represent lattice paths from the origin $(0,0)$ to the point $(p,q)$ on the graph $D_{p,q}$ by a sequence of letters in $\{E,N\}$, e.g., $P = e_1 e_2 \ldots e_{p+q}$ where $e_i \in \{E,N\}$. In particular, any such $P$ must have precisely $p$ $E$-steps and $q$ $N$-steps. Now consider a pair $(\bsy{a},\bsy{b})$ of finite non-negative integer sequences with $\bsy{a} = (a_1, a_2, \ldots, a_p)$ and $\bsy{b} = (b_1, b_2, \ldots, b_q)$. We will say that the order statistics of $(\bsy{a},\bsy{b})$ are bounded by $P$ with respect to the weight set $\bsy{U}$ if and only if, when $r = 1, 2, \ldots, p+q$,
\[ \begin{cases} 
      a_{(i)} < w_U(e_r) & \textnormal{if } e_r \textnormal{ is the $i$-th $E$-step of } P, \\
      b_{(j)} < w_U(e_r)  & \textnormal{if } e_r \textnormal{ is the $j$-th $N$-step of } P,
   \end{cases}
\]
where $a_{(1)} \leq \cdots \leq a_{(p)}$ and $b_{(1)} \leq  \cdots \leq b_{(q)}$ are the order statistics of $\bsy{a}$ and $\bsy{b}$, respectively. 

\begin{definition} \cite{KLY14}  \label{def:U-PF} 
Suppose  the weight set $\bsy{U} =  \{(u_{i,j},v_{i,j}) : 0 \leq i \leq p, 0\leq j\leq q \} \subset \mathbb{N}^2$  satisfies the condition that  $u_{i,j} \leq u_{i',j'}$ and $v_{i,j} \leq v_{i',j'}$ when $(i,j) \preceq (i',j')$. A pair $(\bsy{a},\bsy{b})$ of non-negative integer sequences of respective lengths $p$ and $q$ is a \emph{2-dimensional $\bsy{U}$-parking function} of length $(p,q)$ if and only if the order statistics of $(\bsy{a},\bsy{b})$ are bounded by some (possibly more than one) lattice path from $(0,0)$ to $(p,q)$ with respect to $\bsy{U}$.
\end{definition}

Note that in the above definition, the entries $\{u_{p,j},  v_{i, q}: \ 0 \leq i \leq p,\  0\leq j \leq q\}$ are not used as edge-weight in $D_{p,q}(\bsy{U})$. 
In the following discussion, we will ignore those values and  say that $\bsy{U}$ and $\bsy{U}'$ are the same weight set  whenever they agree  on all the edges of $D_{p,q}$.

\begin{example}
    Figure~\ref{fig:UPF} shows the graph $D_{3,3}(\bsy{U)}$  with weight set 
    $\bsy{U} = \{(u_{i,j},v_{i,j}): 0 \leq i \leq 3, 0 \leq j \leq 3\}$ 
    defined by $u_{i,j} = i+j+1$ and $v_{i,j} = i^2+1$.  
    The pair $(\bsy{a},\bsy{b})$ with $\bsy{a} = (3,0,1)$ and $\bsy{b} = (4,7,4)$ is bounded by the  lattice path in bold  $P = E E N  N E  N$.  Therefore, $(\bsy{a}, \bsy{b})$ is 
    a $\bsy{U}$-parking function of length $(3, 3)$. 
\end{example}

\begin{figure}[H]
\begin{center}
\begin{tikzpicture}[
            > = stealth,
            shorten > = 2pt,
            scale=1.3
        ]
        
        \draw[-] (-.25,-.15) node {\footnotesize{(0,0)}};
        \draw[-] (3.25,3.15) node {\footnotesize{$(3,3)$}};

        \foreach \x in {0,1,2,3} {
        \foreach \y in {0,1,2,3} {
        \node[circle,fill=black,scale=.25] at (\x,\y) {};
        }
        }
        \foreach \x in {0,1,2} {
    \foreach \y in {0,1,2,3} {
       \draw[dotted,->] (\x,\y)--(\x+1,\y); 
      \node[font=\scriptsize] at (\x+.5,\y) {$\the\numexpr\x+\y+1\relax$};
    }
  }

        \foreach \x in {0,1,2,3} {
    \foreach \y in {0,1,2} {
       \draw[dotted,->] (\x,\y)--(\x,\y+1); 
      \node[font=\scriptsize] at (\x,\y+.5) {$\the\numexpr\x*\x+1\relax$};
    }
  }

        \draw[thick,red] (2,0) -- (2,2); 
        \draw[thick,red] (3,2) -- (3,3);
        \draw[thick,red] (0,0) -- (2,0);
        \draw[thick,red] (2,2) -- (3,2);

\end{tikzpicture}
\caption{An example of  $D_{3,3}(\bsy{U})$  and 
a lattice path $P = E E  N N  E N$.} \label{fig:UPF} 
\end{center}
\end{figure}

Let $\cpf_{p,q}^{(2)}(\bsy{U})$  be the set of all 2-dimensional $\bsy{U}$-parking functions with length $(p,q)$.  Although we will only refer to the 2-dimensional case in the present paper, it should be noted that the $k$-dimensional $\bsy{U}$-parking functions for $k > 2$ may be defined in an analogous and straightforward manner. Also, as suggested by the dependency on order statistics in the definition of 2-dimensional $\bsy{U}$-parking functions, the set $\cpf_{p,q}^{(2)}(\bsy{U})$ is invariant under the action of the product of symmetric groups $\mathfrak{S}_p \times \mathfrak{S}_q$. That is, $(\sigma(\bsy{a}),\tau(\bsy{b})) \in \cpf_{p,q}^{(2)}(\bsy{U})$  
for any $(\bsy{a}, \bsy{b}) \in \cpf_{p,q}^{(2)}(\bsy{U})$, $\sigma \in \mathfrak{S}_p$ and $\tau \in \mathfrak{S}_q$. 

Clearly, the set $\cpf_{p,q}^{(2)}(\bsy{U})$ is determined by the set  $\cMPF_{p,q}^{(2)}(\bsy{U})$  of maximal 2-dimensional $\bsy{U}$-parking. 
We say that $(\bsy{a}, \bsy{b}) \in \cpf_{p,q}^{(2)}(\bsy{U})$ is 
\textit{increasing} if both $\bsy{a}$ and $\bsy{b}$ are non-decreasing. 
By permutation invariance, the set $\cMPF_{p,q}^{(2)}(\bsy{U})$ is determined by the set of increasing maximal 2-dimensional 
$\bsy{U}$-parking functions, which is denoted by  $\mathcal{MI}_{p,q}^{(2)}(\bsy{U})$. 
We will often skip the word ``2-dimensional" when there is no confusion.  

\section{The case with affine weight sets} 
\label{sec:affine}

Now that we have recalled all necessary definitions related to the three generalized parking functions that concern this paper, we shift our focus to the overlap between graphical parking functions and 2-dimensional $\bsy{U}$-parking functions.  

For  clarity, in the following of the paper, every graph $G=(V(G),E(G))$ is assumed to be connected with $V(G) = \{0\} \cup A \cup B$, where 0 is the unique root vertex of $G$, $A = \{1,2,\ldots,p\}$, and $B = \{p+1,p+2,\ldots,p+q\}$ for some $p,q \geq 1$. Thus, $|V(G)| = p+q+1$. Every $G$ has an associated edge-weight function $wt_G: E(G) \to \mathbb{Z}^+$. We also write  $wt_G(e) = 0$ to indicate no edge. For $S=A$ or $S=B$, we write $G_S$ to denote the subgraph of $G$ induced by the vertices $\{0\} \cup S$. For any subset $S \subseteq V(G)$, $G[S]$ denotes the subgraph of $G$ induced by $S$.

A $G$-parking function is a sequence of length $p+q$ indexed by the non-zero vertices. For convenience, we will write a $G$-parking function as $(\bsy{a}, \bsy{b})$, where $\bsy{a}$ contains the first $p$ entries  indexed by $A$, and $\bsy{b}$ contains the last $q$ entries indexed by $B$. If $H$ is a subset of $[p+q]$, 
then $\bsy{a}_H$ (resp. $\bsy{b}_H$) is the subsequence of $\bsy{a}$ (resp. $\bsy{b}$) consisting of all the entries indexed by $H \cap A$, 
(resp. $H \cap B$). 

Under the above notation, the main objective of this paper is to answer the following question. 

\noindent 
    \textbf{Question}. For which graphs $G$ and weight sets  $\bsy{U}$ do we have $\cpf_{p,q}^{(2)}(\bsy{U}) = \cpf(G)$?

\medskip 
Note that if $p=0$, then a $\bsy{U}$-parking function has the form $(\emptyset, \bsy{b})$, where $\bsy{b}$ is a $\bsy{v}$-parking function for $\bsy{v}=(v_{0,0}, v_{0,1}, \dots, v_{0, q-1})$.  Hence, it reduces to the case of vector parking functions that is covered by Theorem~\ref{thm:gaydarov and hopkins}. The case $q=0$ is similar. In the following  we assume $p, q>0$. 

In this section, we exam the  special  case where the weight set $\bsy{U}$ is affine, and determine for which $\bsy{U}$ there is a graph $G$ such that  $\cpf_{p,q}^{(2)}(\bsy{U}) = \cpf(G)$. 
The affine case is  particularly significant, as the corresponding Gon\v{c}arov polynomials form a higher dimensional analog of Abel polynomials 
and can be computed efficiently; See \cite{LTY17,LY16}.  
In the special case with $b=d=0$ and $c=c'=a=e=1$ in \Cref{def:affine}, 
the affine  $\bsy{U}$-parking functions coincide with the $(p,q)$-parking functions introduced by Cori and Poulalhon's (\cite{CP02}), as shown by the present authors in \cite[Theorem 3.1]{SY22}. Moreover, 
Cori and Poulalhon proved in 
 \cite[Proposition 7]{CP02}  that $(p,q)$-parking functions are exactly $G$-parking functions for the complete tripartite graph $G=K_{1,p,q}$  with partition sizes $1, p, q$.  As can be seen later,  this result is a special case of Theorem \ref{affinecasetheorem} below. 

 Although the results in this section do not characterize all possible graphs, they offer insights into the underlying structure of such graphs and help us get the complete characterization in Sections \ref{sec:classification} and \ref{sec:finish}. 

\begin{definition}\label{def:affine}
We say that a weight set  $\bsy{U}$ is \emph{affine} if there exists a $2 \times 2$ integer matrix $\bigl( \begin{smallmatrix} b & c \\ c' & d \end{smallmatrix} \bigr)$ 
and integers $a, e$ such that each pair $(u_{i,j},v_{i,j}) \in \bsy{U}$ may be expressed as 
\begin{eqnarray} \label{Eq:matrix}
\begin{pmatrix} 
u_{i,j} \\
v_{i,j}
\end{pmatrix} 
= \begin{pmatrix} 
b & c  \\
c' & d 
\end{pmatrix} 
\begin{pmatrix} 
i \\  j 
\end{pmatrix} + 
\begin{pmatrix} 
a \\
e 
\end{pmatrix}  
\end{eqnarray} 
for $0 \leq i \leq p, 0\leq j \leq q$. 
\end{definition} 
By \Cref{def:U-PF}, if the weight set $\bsy{U}$ is given by \eqref{Eq:matrix}, then  $ b, c, c', d$ are non-negetive. In addition,  at least one of $a$ and $e$ is positive;   otherwise, $u_{0,0}=v_{0,0}=0$ and $\cpf_{p.q}^{(2)}(\bsy{U}) =\emptyset$. We always assume that these conditions are met.

First,  we present a lemma for the  case that $\bsy{U} = \{(u_i,v_j) : 0 \leq i \leq p, 0 \leq j \leq q\} \subset \mathbb{N}^2$, where  $\bsy{u} = (u_0,u_1,\ldots,u_{p-1})$ and $\bsy{v} = (v_0,v_1,\ldots,v_{q-1})$ are vectors of arbitrary  
non-decreasing positive integers.  
The lemma follows directly from 
 Gaydarov and Hopkins' classification of  in Theorem~\ref{thm:gaydarov and hopkins}.  

\begin{lemma}
\label{disjoint-union}
 Let $\bsy{u} = (u_0,u_1,\ldots,u_{p-1})$ and $\bsy{v} = (v_0,v_1,\ldots,v_{q-1})$ each be vectors of one of the forms given by Theorem~\ref{thm:gaydarov and hopkins}(i)-(iii). Let $G_1$ and $G_2$ be graphs such that $\cpf(\bsy{u}) = \cpf(G_1)$ and $\cpf(\bsy{v}) = \cpf(G_2)$. If $\bsy{U} = \{(u_i,v_j) : 0 \leq i < p, 0 \leq j < q \}$, then $\cpf^{(2)}_{p,q}(\bsy{U}) = \cpf(G_1 \cupdot G_2)$, where $G_1 \cupdot G_2$ is the graph obtained by taking the disjoint union of $G_1$ and $G_2$ and merging the two roots into a new root.
\end{lemma}

\begin{proof} It suffices to show that the set $\cMPF^{(2)}_{p,q}(\bsy{U})$ of maximal $\bsy{U}$-parking functions is equal to the set $\cMPF(G)$ of maximal $G$-parking functions, where $G = G_1 \cupdot G_2$. Recall the bijection $\varphi_G: \mathcal{A}(G) \to \cMPF(G)$ given in Proposition~\ref{burningalgorithm}. If $\mathcal{O} \in \mathcal{A}(G)$, the restriction of $\mathcal{O}$ to $G_i$ yields an orientation $\mathcal{O}_i \in \mathcal{A}(G_i)$, for both $i=1, 2$. 
Conversely, any two orientations $\mathcal{O}_1 \in \mathcal{A}(G_1)$ and $\mathcal{O}_2 \in \mathcal{A}(G_2)$ clearly correspond to a unique $\mathcal{O} \in \mathcal{A}(G)$. Let $\psi: \mathcal{A}(G) \to \mathcal{A}(G_1) \times \mathcal{A}(G_2)$ denote this bijection. Since the 2-dimensional $\bsy{U}$-parking functions are precisely the pairs of sequences $(\bsy{a},\bsy{b})$ where $\bsy{a}$ is a $\bsy{u}$-parking function and $\bsy{b}$ is a $\bsy{v}$-parking function, we have the following correspondence:
$$ \cMPF^{(2)}_{p,q}(\bsy{U}) = \cMPF(\bsy{u}) \times \cMPF (\bsy{v}) = \cMPF(G_1) \times \cMPF(G_2)$$ 
$$ \xrightarrow{\varphi_{G_1}^{-1} \times \varphi_{G_2}^{-1}} \mathcal{A}(G_1) \times \mathcal{A}(G_2) \xrightarrow{\psi^{-1}} \mathcal{A}(G) \xrightarrow{\varphi_G} \cMPF(G), $$
where $\cMPF(\bsy{u})$ and $\cMPF(\bsy{v})$ are the sets of maximal $\bsy{u}$-parking functions and maximal $\bsy{v}$-parking functions, respectively. 
\end{proof}

We now provide a characterization of the increasing maximal $\bsy{U}$-parking functions, which play a fundamental role in the discussion to come. 
Each lattice path  $P$ from $(0,0)$ to $(p,q)$ determines an increasing $\bsy{U}$-parking function $(\bsy{a}(P),\bsy{b}(P))$  
as follows: 
if $(u_{0,j_0},u_{1,j_1},\ldots,u_{p-1,j_{p-1}})$ and $(v_{i_0,0},v_{i_1,1},\ldots,v_{i_{q-1},q-1})$ are the sequences of weights on the $E$-steps and the $N$-steps, respectively, for the lattice path $P$, then the corresponding increasing  $\bsy{U}$-parking function is 
\begin{equation} \label{eq:ab_P}
(\bsy{a}(P),\bsy{b}(P)) = (u_{0,j_0}-1,u_{1,j_1}-1,\ldots,u_{p-1,j_{p-1}}-1; v_{i_0,0}-1,v_{i_1,1}-1,\ldots,v_{i_{q-1},q-1}-1).
\end{equation}
By definition, any $\bsy{U}$-parking function is bounded by a pair $(\sigma(\bsy{a}(P)), \tau(\bsy{b}(P)))$ for some lattice path $P$
and permutations $\sigma \in \mathfrak{S}_p$, $\tau \in \mathfrak{S}_q$. It follows that  increasing maximal $\bsy{U}$-parking functions are 
precisely the maximal elements  in the set $\{(\bsy{a}(P), \bsy{b}(P)):  P \text{ is a lattice path from $(0,0)$ to  $(p,q)$}\}$.

To describe the overlap between graphical parking functions and 2-dimensional $\bsy{U}$-parking functions for affine $\bsy{U}$, 
we distinguish three cases, according to the values of $c$ and $c'$. 

\medskip 

\noindent \textbf{Case 1. $c=c'=0$. } \ 
Then $u_{i,j}=a+bi$  and $v_{i,j}=e+dj$.   If $a=0$, then $u_{0,j}=0$ for all $j$, which implies $\cpf_{p,q}^{(2)}(\bsy{U})=\emptyset$, and there is no corresponding graph. A similar argument holds for  $e=0$. If $a \geq 1$ and $e \geq 1$, then there exists a graph $G$ with $\cpf(G)= \cpf_{p,q}^{(2)}(\bsy{U})$, as described in \Cref{disjoint-union}.   

\medskip 
\noindent \textbf{Case 2.} Only one of $c, c'$ is 0. \ 
Without loss of generality, assume  $c=0$ but $c' \geq 1$. 
Let $P_1=E^pN^q$. We observe that for any lattice path $P$ from $(0,0)$ to $(p,q)$, $(\bsy{a}(P), \bsy{b}(P)) \preceq (\bsy{a}(P_1), \bsy{b}(P_1))$. Hence there is only one increasing maximal 
$\bsy{U}$-parking functions. Let 
$\bsy{U'} = \{(u'_{i,j},v'_{i,j}): 0 \leq i \leq p,\  0 \leq j \leq q\} $ is given by 
\begin{eqnarray} \label{Eq:matrix'}
\begin{pmatrix} 
u'_{i,j} \\
v'_{i,j}
\end{pmatrix} 
= \begin{pmatrix} 
b & 0  \\
0 & d 
\end{pmatrix} 
\begin{pmatrix} 
i \\  j 
\end{pmatrix} + 
\begin{pmatrix} 
a \\
e+c'p 
\end{pmatrix}.  
\end{eqnarray} 
Then  $(\bsy{a}(P_1), \bsy{b}(P_1))$   is also the only increasing maximal $\bsy{U}'$-parking function.  
Consequently, we have $\cMPF_{p,q}^{(2)}(\bsy{U}) = \cMPF_{p,q}^{(2)}(\bsy{U'})$,  and hence,  $\cpf_{p,q}^{(2)}(\bsy{U}) = \cpf_{p,q}^{(2)}(\bsy{U'})$. 
Now the problem is reduced to Case 1,  and 
there exists a graph $G$ with $\cpf(G)= \cpf_{p,q}^{(2)}(\bsy{U})$ if and only if  $a \geq 1$. 

By symmetry, when $c \geq 1$ and $c'=0$, 
there exists  a graph $G$ with $\cpf(G)= \cpf_{p,q}^{(2)}(\bsy{U})$ if and only if  $e \geq 1$. 

\medskip 
\noindent \textbf{Case 3.} $c \geq 1$ and $c' \geq 1$.   
First we prove a lemma. 

\begin{lemma} \label{lemma:c=c'}
     Suppose $\bsy{U} = \{(u_{i,j},v_{i,j}): 0 \leq i \leq p,\  0 \leq j \leq q\}$ is  given by \eqref{Eq:matrix} with $c\geq 1$ 
    and $c' \geq 1$. If there exists a graph $G$ such that 
    $\cpf_{p,q}^{(2)}(\bsy{U})  = \cpf(G)$, then $c=c'$. 
\end{lemma}
\begin{proof}
     In the weighted digraph $D_{p,q}(\bsy{U})$, the horizontal edge starting at $(i,j)$  has weight $u_{i,j}=bi+cj+a$, while the vertical edges starting at $(i,j)$ has weight $v_{i,j}=c'i+dj+e$.  
    When both $c\geq 1$ and $c' \geq 1$, we have $u_{i,j} < u_{i, j'}$ whenver $j < j'$, and $v_{i,j} < v_{i',j}$ whenever $j < j'$.  
    Therefore   
    $(\bsy{a}(P), \bsy{b}(P))$ and $(\bsy{a}(Q), \bsy{b}(Q))$ are non-comparable under $\preceq$ for any two distinct lattice paths $P, Q$ from $(0,0)$ to $(p,q)$. 
     Hence the map $P \rightarrow (\bsy{a}(P), \bsy{b}(P))$ 
      is a bijection between lattice paths from $(0,0)$ to $(p,q)$ in $D_{p,q}$ and 
       increasing maximal $\bsy{U}$-parking functions. 
      


For the two lattice paths $P_1=E^pN^q$ and $P_2=N^q E^p$, the increasing maximal $\bsy{U}$-parking functions 
$(\bsy{a}(P_1), \bsy{b}(P_1)$  
and 
$(\bsy{a}(P_2), \bsy{b}(P_2))$ are given by 
\begin{eqnarray*}
    \bsy{a}(P_1) &= &(a-1, a+b-1, \cdots, a+b(p-1) -1), \\ 
    \bsy{b}(P_1) &= & (c'p+e-1, c'p+d+e-1, \cdots, c'p+d(q-1)+e-1), \\
    \bsy{a}(P_2) & = & (a+cq-1, a+b+cq-1, \cdots, a+b(p-1)+cq-1),   \\ 
    \bsy{b}(P_2)  & = & (e-1, e+d-1, \cdots, e+d(q-1)-1).     
\end{eqnarray*}

An immediate consequence of Proposition \ref{burningalgorithm} 
     is that the sum of the entries of any maximal $G$-parking function is always the same. This is because for any orientation in $\mathcal{A}(G)$, the sum of all indegrees is exactly the total weight on all the edges.  Hence if 
    $\cpf_{p,q}^{(2)}(\bsy{U})  = \cpf(G)$ for a graph $G$, then 
    the  sum of all entries of $\bsy{a}(P_1)$ and $\bsy{b}(P_1)$ should be equal to the sum of entries of $\bsy{a}(P_2)$ and $\bsy{b}(P_2)$. This leads to $c=c'$.
\end{proof}

Thus, if $ c,c' \geq 1$ and $c\neq c'$ in \eqref{Eq:matrix}, then there is no graph $G$ such that $\cpf_{p,q}^{(2)} (\bsy{U}) =\cpf(G)$.  
The case $c=c' \geq 1$ is described by the next theorem. 

\begin{theorem}
\label{affinecasetheorem}
 Suppose $\bsy{U} = \{(u_{i,j},v_{i,j}): 0 \leq i \leq p, 0 \leq j \leq q\} \subset \mathbb{N}^2$ is  given by 
  $(u_{i,j},v_{i,j})^{\textnormal{T}} = \bigl( \begin{smallmatrix} b & c \\ c & d \end{smallmatrix} \bigr) \bigl( \begin{smallmatrix} i \\ j \end{smallmatrix} \bigr) + \bigl( \begin{smallmatrix} a\\ e\end{smallmatrix} \bigr)$ with $c \in \mathbb{Z}^+$,  $a, b, d, e \in \mathbb{N}$,  and at least one of $a, e$ is not $0$, then $\cpf^{(2)}_{p,q}(\bsy{U}) = \cpf(G)$ where $G = K_{p+q+1}$ with vertex set $[p+q]_0$ and edge-weight function
\[ wt_G(\{i,j\}) = \begin{cases} 
	  a & \textnormal{if } i = 0 \textnormal{ and } j = 1,\ldots,p; \\
	  b & \textnormal{if } 1 \leq i < j \leq p; \\
	  c &  \textnormal{if } 1 \leq i \leq p \textnormal{ and } p+1 \leq j \leq p+q; \\
      d & \textnormal{if } p+1 \leq i < j \leq p+q; \\
      e & \textnormal{if } i = 0 \textnormal{ and } j = p+1,\ldots,p+q.
   \end{cases} 
\]
 
\end{theorem}

\begin{proof}

We will show that the map $\varphi: \mathcal{A}(G) \to \cMPF^{(2)}_{p,q}(\bsy{U})$ given by
$$\mathcal{O} \mapsto (\indeg_{_\mathcal{O}}(1)-1, \ldots, \indeg_{_\mathcal{O}}(p)-1; \indeg_{_\mathcal{O}}(p+1)-1, \ldots, \indeg_{_\mathcal{O}}(p+q)-1)$$
is a bijection. By Proposition \ref{burningalgorithm}, this  implies $\cMPF(G) = \cMPF^{(2)}_{p,q}(\bsy{U})$, and we are be done.

Fix $\mathcal{O} \in \mathcal{A}(G)$, and let $S = \{1,\ldots,p\}$ and $T = \{p+1,\ldots,p+q\}$. From $\mathcal{O}$, we construct a lattice path $P$ from $(0,0)$ to $(p,q)$ by the following process:
\begin{enumerate}
	\item Delete vertex 0 and all its incident edges from $G$.
	\item Choose a vertex $v$ with $\indeg_{\mathcal{O}}(v)=0$. If more than one such vertex exists, choose $v$ among these at random.  (It does not affect the remaining steps.)     
    If $v \in S$, let the next step of $P$ be an $E$-step, and if $v \in T$, let the next step of $P$ be an $N$-step.
	\item Delete $v$ and all its incident edges. If a non-empty subgraph remains, return to step 2.
\end{enumerate}
Note that since $\mathcal{O}$ is acyclic, there will always exist a vertex $v$ having $\indeg_{\mathcal{O}}(v)=0$ in step 2. Now we claim that the order statistics of the sequence 
$$(\indeg_{_\mathcal{O}}(1)-1, \ldots, \indeg_{_\mathcal{O}}(p)-1; \indeg_{_\mathcal{O}}(p+1)-1, \ldots, \indeg_{_\mathcal{O}}(p+q)-1)$$ 
are bounded by $P$ with respect to the weight $\bsy{U}$. For if $v \in S$ corresponds to a step from $(i,j)$ to $(i+1,j)$ in $P$, this means there were $i$ vertices from $S$ removed from $G$ prior to the removal of $v$ according to the above process, and every edge between one of these $i$ vertices and $v$ was directed towards $v$. Additionally, there were $j$ vertices from $T$ removed from $G$ prior to the removal of $v$, and every edge from each of these $j$ vertices to $v$ was directed towards $v$. Hence, $\indeg_{_\mathcal{O}}(v) - 1 = a+bi + cj - 1 < a+bi + cj = u_{i,j}$. Similarly, if $v \in T$ corresponds to a step from $(i,j)$ to $(i,j+1)$ in $P$, this same argument implies that $\indeg_{_\mathcal{O}}(v) - 1 = e+ci + dj - 1 < e+ci + dj = v_{i,j}$. This proves that $(\indeg_{_\mathcal{O}}(1)-1, \ldots, \indeg_{_\mathcal{O}}(p)-1; \indeg_{_\mathcal{O}}(p+1)-1, \ldots, \indeg_{_\mathcal{O}}(p+q)-1)$ is indeed a $\bsy{U}$-parking function. Moreover, it is maximal as 
its order statistics are  exactly $(\bsy{a}(P), \bsy{b}(P))$, which is an increasing maximal $\bsy{U}$-parking function. 
Hence $\varphi$ is a well-defined map. 

To show injectivity, note that every orientation $\mathcal{O} \in \mathcal{A}(G)$ is uniquely determined by its indegree sequence: starting with the sink vertices (i.e., those $v$ for which $\indeg_{\mathcal{O}}(v) = deg(v)$), orient all their adjacent vertices towards the sink, then remove the sink vertices, yielding a subgraph of $G$. Repeat this process until the resulting subgraph of $G$ only consists of the root 0. Thus, $\varphi(\mathcal{O}_1) \neq \varphi(\mathcal{O}_2)$ for any two distinct $\mathcal{O}_1,\mathcal{O}_2 \in \mathcal{A}(G)$. 

To show surjectivity, fix a sequence $\bsy{x} = (a_1,\ldots,a_p ; a_{p+1},\ldots,a_{p+q}) \in \cMPF^{(2)}_{p,q}(\bsy{U})$. 
Let $ b_1< \cdots <  b_p$ be the order statistics of $a_1, \dots, a_p$ and 
$c_1 < \cdots < c_q$ be the order statistics of $a_{p+1}, \dots, a_{p+q}$. 
Let $\sigma \in \mathfrak{S}_p$ be determined by $b_i=a_{\sigma(i)}$, and similarly, $\tau \in \mathfrak{S}_q$ be determined by $c_j=a_{p+\tau(j)}$. 
Let $\bsy{y}=(b_1, \dots, b_p; c_1, \dots, c_q)$, which is an increasing maximal $\bsy{U}$-parking functions. 
Then $\bsy{y}= (\bsy{a}(P), \bsy{b}(P))$ for a unique lattice path 
$P = e_1 e_2 \ldots e_{p+q}$ from $(0,0)$ to $(p,q)$, in which  $e_r \in \{E,N\}$ for each $r$.  

The following process assigns an orientation $\mathcal{O}$ to $G=K_{1+p+q}$. 
\begin{enumerate}
	\item Direct all edges incident to 0 away from 0.
	\item If $e_1 = E$, direct all non-oriented edges incident to vertex $\sigma(1)$ away from $\sigma(1)$; if $e_1 = N$, direct all non-oriented edges incident to vertex $p+\tau(1)$ away from $p+\tau(1)$.
	\item Consider the next step $e_r$ of $P$. If $e_r$ is an $E$-step from $(i,j)$ to $(i+1,j)$, direct all non-oriented edges incident to vertex $\sigma(i+1)$ away from $\sigma(i+1)$; if $e_r$ is an $N$-step from $(i,j)$ to $(i,j+1)$, direct all non-oriented edges incident to vertex $p+\tau(j+1)$ away from $p+\tau(j+1)$.
	\item If there remain non-oriented edges of $G$, return to step 3.
\end{enumerate}
It is clear that this process assigns a valid orientation $\mathcal{O}$ to $G$. To see that $\mathcal{O}$ is acyclic, suppose there existed a cycle $v_1, \ldots, v_k, v_1$. Then step 3 of the process above would imply that the vertices $v_1, \ldots, v_k$ are encountered in that order. We obtain a contradiction, though, if there is a directed edge from $v_k$ to $v_1$. It is also clear that 0 is the unique vertex with $\indeg_{_\mathcal{O}}(0)=0$. Thus, $\mathcal{O} \in \mathcal{A}(G)$, and so $\varphi$ is a bijection. \end{proof}

To illustrate the process of going from an orientation $\mathcal{O} \in \mathcal{A}(G)$ to a lattice path $P$ described in the above proof, consider the following example.

\begin{example}
Let $p=3$, $q=2$, and $G = K_{1+p+q}$ with edge-weights indicated in Figure \ref{fig:lattice_path_example} below.  For simplicity, we assume $b=d=0$.
Let $\mathcal{O} \in \mathcal{A}(G)$ be the acyclic orientation shown on the left. The corresponding lattice path $P=ENENE$ is given in bold on the right, and the maximal $\bsy{U}$-parking function $\varphi(\mathcal{O})$ is $((a+c-1, a+2c-1, a-1), (2c+e-1, c+e-1))$.  For this maximal $\bsy{U}$-parking function, $\sigma=312$ and $\tau=21$. 
\end{example}

\begin{figure}[H]
\begin{center}
\begin{subfigure}[t]{.5\textwidth}
\begin{center}
\begin{tikzpicture}[line width=.7pt,scale=1.1]
\tikzstyle{vertex}=[draw,circle,fill=black,minimum size=3pt,inner sep=0pt]
\path (0,0) node[vertex, label=left:0] (v0) {};
\path (1.25,1) node[vertex, label=above:1] (v1) {};
\path (2.5,1) node[vertex, label=above:2] (v2) {};
\path (3.75,1) node[vertex, label=above:3] (v3) {};
\path (1.875,-1) node[vertex, label=below:4] (v4) {};
\path (3.125,-1) node[vertex, label=below:5] (v5) {};
\draw[-Stealth, blue] (v0) -> (v1) node [midway, above=0.5mm] {$a$};
\draw[-Stealth, blue] (v0) -> (v2);
\draw[-Stealth, blue] (v0) -> (v3);
\draw[-Stealth, Green] (v0) -> (v4) node [midway, below=0.5mm] {$e$};
\draw[-Stealth, Green] (v0) -> (v5);
\draw[-Stealth, red] (v1) -> (v4) node [midway, left=0.5mm] {$c$};
\draw[-Stealth, red] (v5) -> (v1);
\draw[-Stealth, red] (v4) -> (v2);
\draw[-Stealth, red] (v5) -> (v2);
\draw[-Stealth, red] (v3) -> (v4);
\draw[-Stealth, red] (v3) -> (v5);
\end{tikzpicture}
\caption*{}
\label{fig:sub1}
\end{center}
\end{subfigure}%
\begin{subfigure}[t]{.5\textwidth}
\begin{center}
\begin{tikzpicture}[line width=.7pt,scale=1.6]
\draw[step=1cm,dotted,color=gray] (0,0) grid (3,2);
\path (0.5,0.15) node {$a$};
\path (1.5,0.15) node {$a$};
\path (2.5,0.15) node {$a$};
\path (0.5,1.15) node {$c+a$};
\path (1.5,1.15) node {$c+a$};
\path (2.5,1.15) node {$c+a$};
\path (0.5,2.15) node {$2c+a$};
\path (1.5,2.15) node {$2c+a$};
\path (2.5,2.15) node {$2c+a$};
\path (-0.15,0.5) node {$e$};
\path (-0.15,1.5) node {$e$};
\path (0.9,0.5) node {$c+e$};
\path (0.9,1.5) node {$c+e$};
\path (1.88,0.5) node {$2c+e$};
\path (1.88,1.5) node {$2c+e$};
\path (2.88,0.5) node {$3c+e$};
\path (2.88,1.5) node {$3c+e$};
\draw [ultra thick, black,scale=1] (0,0) -- (1,0) -- (1,1) -- (2,1) -- (2,2) -- (3,2);
\draw [<->,gray, scale=1] (0,2.5) -- (0,0) -- (3.5,0);
\end{tikzpicture}
\caption*{}
\label{fig:sub2}
\end{center}
\end{subfigure}
\end{center}
\caption{An $\mathcal{O} \in \mathcal{A}(G)$  and its corresponding lattice path $P$.}
\label{fig:lattice_path_example}
\end{figure}
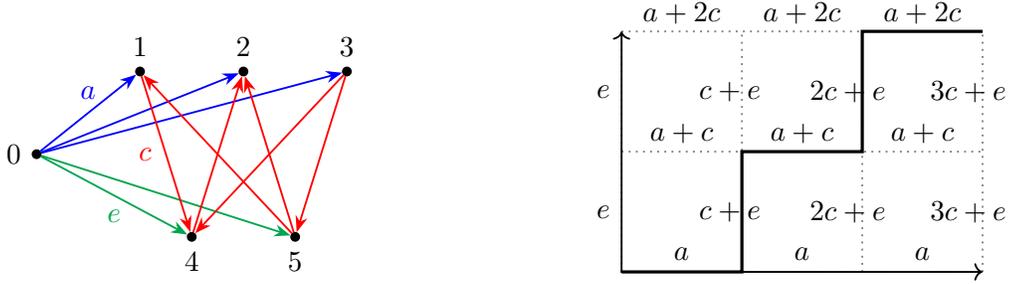

\section{Classification of $G$ for which $\cpf(G)$ is $(\mathfrak{S}_p \times \mathfrak{S}_q)$-invariant}
\label{sec:classification}

In this section, we present a total classification of all graphs $G$ whose set of $G$-parking functions $\cpf(G)$ is invariant under the action of the product of symmetric groups $\mathfrak{S}_p \times \mathfrak{S}_q$. Clearly, if there exists a 2-dimensional set of weights  $\bsy{U}$ for which $\cpf_{p,q}^{(2)}(\bsy{U}) = \cpf(G)$, a necessary condition is that any $\mathfrak{S}_p \times \mathfrak{S}_q$-permutation of a $G$-parking function is itself a $G$-parking function. In fact, we shall see in \Cref{sec:finish} that every graph satisfying this invariance property has the stated equivalence for some $\bsy{U}$.

Before stating the main result of this section, we define some specialized graphs. We will say that an $a$-tree $G$ is an \textit{$a$-star} if $G$ is the complete bipartite graph $K_{1,n}$, where $n = |V(G)| - 1$ and 0 is the internal vertex. If the $a$-tree $G$ is a path, we will say $G$ is an \textit{$a$-path}. Lastly, recall that if $G=T$ is a tree, then there is exactly one (acyclic) orientation of $T$ in which 0 is the unique source vertex, so that $|\mathcal{A}(T)| = 1$. Letting $\mathcal{A}(T) = \{\mathcal{O}\}$, we can define $wt_T: E(G) \to \mathbb{Z}^+$ to satisfy
\begin{equation} \label{Eq:Tweight}
wt_T(\{i,j\}) = \begin{cases} 
      a & \textnormal{if } (i,j) \in \mathcal{O} \textnormal{ and } j \in A \\
     b & \textnormal{if } (i,j) \in \mathcal{O} \textnormal{ and } j \in B 
   \end{cases} 
\end{equation}
for any $\{i,j\} \in E(T)$ and $a,b \in \mathbb{Z}^+$. Then we will say the tree $T$ is an \textit{$(a,b)$-tree} when $T$ is equipped with the edge-weight function $wt_T$ given in \eqref{Eq:Tweight}.

 The  main result is stated in the following theorem. 

\begin{theorem}
\label{maintheorem}
Fix positive integers $p$ and $q$. 
Let $G$ be a graph with $V(G) = \{0\} \cup A \cup B$, where $0$ is the unique root vertex of $G$, $A= \{1,2, \dots, p\}$,  and $B = \{p+1, \dots, p+q\}$.  
Let $wt_G: E(G) \to \mathbb{N}$ be the associated edge-weight function of $G$, and assume $N_G(0) \cap A \neq \emptyset$.
Denote by $G_A$ the subgraph of $G$ induced by the vertices $\{0\} \cup A$. 
If $\cpf(G)$ is invariant under the action of $\mathfrak{S}_p \times \mathfrak{S}_q$, then one of the following holds:
\begin{enumerate}
    \item[(i)] $G$ is a cycle, and for $a,b \in \mathbb{Z}^+$ with $a \neq b$:
    \begin{enumerate}
        \item[(a)] $G$ is an $a$-cycle; or
        \item[(b)] $p=1$. $A=\{1\}$ with $wt_G(\{0,1\}) = a$, and $wt_G(\{i,j\}) = b$ for all other $\{i,j\} \in E(G)$; or
        \item[(c)] $p=2$. $A = \{1,2\}$ with $wt_G(\{0,1\}) = wt_G(\{0,2\}) = a$ and $wt_G(\{i,j\}) = b$ for all other $\{i,j\} \in E(G)$.
    \end{enumerate}
    
    \item[(ii)] $p=2$. $A = \{1,2\}$ and $G$ is a cycle with the chord $\{1,2\} \in E(G)$, characterized by the following edge-weight function $(a,b,c \in \mathbb{Z}^+)$:
    \[ wt_G(\{i,j\}) = \begin{cases} 
	  a & \textnormal{if } i = 0 \textnormal{ and } j \in A, \\
	  b & \textnormal{if } i,j \in A, \\
	  c & \textnormal{ otherwise. }  \end{cases}
    \]
    
    \item[(iii)] $G$ is obtained from the complete graph $K_{p+q+1}$ by assigning the edge-weight function below, where $a,c \in \mathbb{Z}^+$ and $b,d,e \in \mathbb{N}$:
    \[ wt_G(\{i,j\}) = \begin{cases} 
	  a & \textnormal{if } i = 0 \textnormal{ and } j \in A, \\
	  b & \textnormal{if } i,j \in A, \\
	  c & \textnormal{if } i \in A \textnormal{ and } j \in B, \\
      d & \textnormal{if } i,j \in B, \\
      e & \textnormal{if } i=0 \textnormal{ and } j \in B.
    \end{cases}
    \]

    \item[(iv)] $G$ is an $(a,b)$-tree for some $a,b \in \mathbb{Z}^+$.
    
    \item[(v)] $G_A$ is an $a$-cycle or $K_{p+1}^{a,b}$ for some $a,b \in \mathbb{Z}^+$, and for $c,d \in \mathbb{Z}^+$:
    \begin{enumerate}
        \item[(a)] there is a vertex $i \in A \cup \{0\}$ such that $G[\{i\} \cup B]$ is a $c$-tree or a $c$-cycle or $K_{q+1}^{c,d}$ (with root $i$), and no other vertices of $A \cup \{0\}$ are adjacent to vertices of $B$; or
        \item[(b)] there are vertices $i_1, i_2, \ldots, i_n \in A \cup \{0\}$ ($n > 1$) such that $G[\bigcup_{k=1}^n \{i_k\} \cup B]$ is a forest of $c$-trees, and no other vertices of $A \cup \{0\}$ are adjacent to vertices of $B$.
    \end{enumerate}

    \item[(vi)] $G_A$ is a forest of $a$-trees for some $a \in \mathbb{Z}^+$, and there is a vertex $i$ in the $0$-component of $G_A$ such that $G[\{i\} \cup B]$ is a $c$-cycle or $K_{q+1}^{c,d}$ (with root $i$) for some $c,d \in \mathbb{Z}^+$, with no cycle of $G$ intersecting $(A \cup \{0\}) \setminus \{i\}$. 
     The none-0 tree component of $G_A$ are rooted at vertices in $B$.
     (See Figure~\ref{fig:case-6} for an illustration. ) 
\end{enumerate}
\end{theorem}

\begin{figure}[H]
\begin{center}
\begin{tikzpicture}[line width=.6pt,scale=.8]
\tikzstyle{vertex}=[draw,circle,fill=black,minimum size=4pt,inner sep=0pt]
\tikzstyle{evertex}=[draw,circle,fill=none,minimum size=5pt,inner sep=2pt,thick]

\path (-.5,0) node[vertex, label=left:0] (v0) {};
\path (1,-1) node[vertex,red] (a1) {};
\path (2,-1.5) node[vertex,red] (a2) {};
\path (3,-1.5) node[vertex,red] (a3) {};
\path (2,-0.5) node[vertex,red] (a4) {};
\path (1,1) node[vertex,red] (a5) {};
\path (2,1) node[vertex,red,label=above:$u$] (a6) {};
\path (3,.5) node[evertex,draw=blue] (b1) {};
\path (4,1) node[evertex, draw=blue] (b2) {};
\path (3,1.5) node[evertex,draw=blue] (b3) {};
\path (4,0) node[vertex,red] (a7) {};
\path (4,2.) node[vertex,red] (a8) {};
\path (5,2) node[vertex,red] (a9) {};

\draw[red] (v0) -- (a1);
\draw[red] (v0) -- (a5) node[red,midway,above] (a) {$a$};
\draw[red] (a6) -- (a5);
\draw[red] (v0) -- (a1) -- (a2) -- (a3);
\draw[red] (a1) -- (a4);
\draw[blue] (a6) edge (b1);
\draw[blue] (a6) -- (b2);
\draw[blue] (a6) -- (b3) node[blue,midway,below=3mm] (b) {$b$};
\draw[Green] (b3) -- (b2) node[Green,midway,above=1.5mm,right=-2mm] (c) {$c$};
\draw[Green] (b2) -- (b1) -- (b3);
\draw[red] (b1) edge (a7);
\draw[red] (b3) -- (a8) -- (a9);

\end{tikzpicture}
\end{center}
\caption{An example of a graph $G$ described in Theorem~\ref{maintheorem}(vi). Red nodes are vertices in $A$, while blue circle-nodes are vertices in $B$. }
\label{fig:case-6}
\end{figure}

Note that it is possible for a graph $G$ to satisfy more than one of the characterizations in Theorem \ref{maintheorem}.  Graphs of the form $G_1 \cupdot G_2$, where $G_1$ and $ G_2$ are graphs described in Section 2.2, are special cases spread in cases (iii)--(vi).

Our approach for proving Theorem~\ref{maintheorem} will be piecemeal: after proving some initial results concerning graphs whose $G$-parking functions have the required invariance property in Subsection \ref{subsec:property}, we  separately consider the case where $G_A$ is connected (with $G_B$ either connected or not) in Subsection \ref{subsec:connect} and the case where neither $G_A$ nor $G_B$ is connected in Subsection \ref{subsec:not-connect}. Our proofs heavily utilize Gaydarov and Hopkins' results on the one dimensional case to characterize the possible structures of both $G_A$ and $G_B$. 

\subsection{Some properties of $G$ and $\cpf(G)$}
\label{subsec:property}

In this subsection, we establish some preliminary results concerning graphs $G$ for which  $\cpf(G)$ is invariant under the action of $\mathfrak{S}_p \times \mathfrak{S}_q$. As in \cite{GH16}, $\mathcal{I}(G) \subseteq V(G)$ will denote the set of vertices in $A \cup B$ which are not cut vertices. Note that since $G$ is assumed to be connected, the deletion of a cut vertex always yields a disconnected graph.

\begin{lemma}
\label{lemma:invariance_maximal_GPF}
For any graph $G$, $\cpf(G)$ is invariant under the action of $\mathfrak{S}_p \times \mathfrak{S}_q$ if and only if $\cMPF(G)$ is invariant under the action of the same group.
\end{lemma}

\begin{proof}  Fix $\sigma, \tau \in \mathfrak{S}_p \times \mathfrak{S}_q$. Proposition~\ref{burningalgorithm} implies that all maximal $G$-parking functions have the same sum of their entries. For $\sigma, \tau \in \mathfrak{S}_p \times \mathfrak{S}_q$, $(\bsy{a},\bsy{b}) \in \cMPF(G)$ and $(\sigma(\bsy{a}), \tau(\bsy{b}))$ have equal sums. Thus, since $(\sigma(\bsy{a}), \tau(\bsy{b})) \in \cpf(G)$, we must have $(\sigma(\bsy{a}), \tau(\bsy{b})) \in \cMPF(G)$. 

Conversely, if $(\bsy{a},\bsy{b}) \in \cpf(G)$, there exists $(\bsy{a}',\bsy{b}') \in \cMPF(G)$ with $\bsy{a} \preceq \bsy{a}'$ and $\bsy{b} \preceq \bsy{b}'$. Since $(\sigma(\bsy{a}'),\tau(\bsy{b}')) \in \cMPF(G)$, it follows that $(\sigma(\bsy{a}), \tau(\bsy{b})) \in \cpf(G)$. 
\end{proof} 

\begin{lemma}
\label{lemma:non_cut_vertex_equal_degrees}
Let $G$ be a graph for which $\cpf(G)$ is invariant under the action of $\mathfrak{S}_p \times \mathfrak{S}_q$. Then for $S=A$ or $S=B$, $\deg_G(i) = \deg_G(i')$ for all $i,i'\in \mathcal{I}(G) \cap S$.
\end{lemma}

\begin{proof} 
 By symmetry, it is sufficient to just prove the statement for $S=A$. Choose $i \in \mathcal{I}(G) \cap A$, and let $G \setminus \{i\}$ be the graph obtained from $G$ by deleting vertex $i$. Since $i$ is not a cut vertex of $G$, $G \setminus \{i\}$ is still connected, and hence $\mathcal{A}(G \setminus \{i\})$ is non-empty. Fix $\mathcal{O}_i \in \mathcal{A}(G \setminus \{i\})$ and extend $\mathcal{O}_i$ into an orientation $\mathcal{O}_{i,G} \in \mathcal{A}(G)$ by directing all edges of $E(G)$ adjacent to $i$ toward $i$. Then $i$ is a sink under the orientation $\mathcal{O}_{i,G}$, with $\indeg_{\mathcal{O}_{i,G}}(i) = \deg_G(i)$. Clearly, $\max_{\mathcal{O} \in \mathcal{A}(G)} \indeg_{\mathcal{O}}(i) = \indeg_{\mathcal{O}_{i,G}}(i) = \deg_G(i)$, which holds for all $i \in \mathcal{I}(G) \cap A$.

Suppose  $i, i'$ are two distinct vertices in $\mathcal{I}(G) \cap A$.  The group $\mathfrak{S}_p$ acts transitively on the first $p$ entries of sequences in $\cMPF(G)$. 
Since $\cMPF(G)$ is invariant under the action of $\mathfrak{S}_p \times \mathfrak{S}_q$, 
we have 
$\max_{\mathcal{O} \in \mathcal{A}(G)} \indeg_{\mathcal{O}}(i) = 
\max_{\mathcal{O} \in \mathcal{A}(G)} \indeg_{\mathcal{O}}(i')$, 
which implies $\deg_G(i)=\deg_G(i').$
\end{proof} 

Lemma~\ref{lemma:non_cut_vertex_equal_degrees} states that any two non-cut vertices which are both in $A$ or both in $B$ in a graph  have the same maximal indegree (equal to their degree) over all orientations in $\mathcal{A}(G)$. This fact will be especially useful when we work with cut-vertex-free graphs.

The next result is a crucial tool in our classification efforts. 
Generally speaking,  if $G_1$ is a connected induced subgraph of a graph $G$ for which $\cpf(G)$ possesses the specified invariance property, then $\cpf(G_1)$ inherits that invariance property as well.

Let $H \neq \emptyset $ be a proper subset of $[p+q]$ and $H^C$ denote the complement of $H$ in $[p+q]$. Assume that $G_H$, the induced subgraph of $G$ on $H \cup \{0\}$, is connected. 
Let $\mathcal{A}(G)_{H \to H^C}$ represent the subset of $\mathcal{A}(G)$ consisting of all orientations in which every edge $\{i,j\} \in E(G)$ between some $i \in H$ and $j \in H^C$ is oriented towards $j$. 
Then  $\mathcal{A}(G)_{H \to H^C}$ is nonempty. 
Specifically, we can take any $\mathcal{O}_{H} \in \mathcal{A}(G_H)$ and extend it into an orientation $\mathcal{O}_G$ of $G$ by 
the following steps: 
\begin{enumerate} 
\item  Directing all edges $\{i,j\}$, where $i \in H$ and $j \in H^C$, towards $j$. 
\item  Order the  vertices in $H^C$ by a breadth search order, starting from the neighbors of $H$. That is, assume that $v_1, \dots, v_t \in H^C$ are neighbors of vertices in $H$. starting from this list of length $t$, for $i=1$ to $|H^c|$, if a vertex $x$ is a neighbor of $v_i$ but not in the list, add $x$ to the end of the list.   This will create an ordering of all vertices in $H^c$. 
Direct an edge $\{i,j\}$ toward $j$ if $j$ appears after $i$ in this ordering. 
\end{enumerate} 
The obtained $\mathcal{O}(G)$ is in $ \mathcal{A}(G)_{H \to H^C}$.

Let $k = |H \cap A|$ and $ l= |H \cap B|$. Then $\cpf(G_H)$ 
can be written as $(\bsy{c}, \bsy{d})$, where $\bsy{c}$ is of length $k$ and $\bsy{d}$ of length $l$.

\begin{proposition}
\label{prop:subgraph_invariance}
Let $G$ be a graph for which $\cpf(G)$ is invariant under the action of $\mathfrak{S}_p \times \mathfrak{S}_q$. Assume $H \subseteq [p+q]$ is a subset of vertices such that $G_H$ is a connected graph with root 0. 
Then 
$\cpf(G_H)$ is invariant under the action of the product of symmetric groups $\mathfrak{S}_k \times \mathfrak{S}_l$.  
\end{proposition}

\begin{proof}
 Given a maximal $G_H$-parking function $(\bsy{c},\bsy{d}) \in \cMPF(G_H)$, let $\mathcal{O}_H \in \mathcal{A}(G_H)$ be the acyclic
 orientation of $G_H$ that corresponds to $(\bsy{c},\bsy{d})$. The argument in the proceeding paragraph shows that 
  there exists an orientation $\mathcal{O} \in \mathcal{A}(G)_{H \to H^C}$ whose restriction to $G_H$ is $\mathcal{O}_H$. 
  Let $(\bsy{a},\bsy{b}) \in \cMPF(G)$ be the maximal $G$-parking function corresponding to $\mathcal{O}$. 
  Then  
 $\bsy{c} = \bsy{a}_H$ and its complement  in $\bsy{a}$ is $\bsy{c}'=\bsy{a}_{H^C}$. Similarly,  $\bsy{d}=\bsy{b}_H$ and its complement in $\bsy{b}$ is $\bsy{d}'=\bsy{b}_{H^C}$.

   Now fix $(\sigma,\tau) \in \mathfrak{S}_k \times \mathfrak{S}_l$.   Extend $\sigma$ to $\sigma' \in \mathfrak{S}_p$ by letting $\sigma'(i)= \sigma(i)$ for $i \in H \cap A$, and $\sigma'(i)=i$ for $i \in H^C \cap A$. Similarly extend $\tau$ to $\tau' \in \mathfrak{S}_q$. 
   By Lemma~\ref{lemma:invariance_maximal_GPF}, $(\sigma'(\bsy{a}),\tau'(\bsy{b})) \in \cMPF(G)$. 
   Let  $\mathcal{O}' \in \mathcal{A}(G)$ be the corresponding orientation for $(\sigma'(\bsy{a}),\tau'(\bsy{b}))$. 
   We now wish to show that $\mathcal{O}' \in \mathcal{A}(G)_{H \to H^C}$. If this is the case, then the restriction of $\mathcal{O}'$ to $G_H$ is an acyclic  orientation of $\mathcal{A}(G_H)$ which gives $(\sigma(\bsy{c}), \tau(\bsy{d})$, and hence, $(\sigma(\bsy{c}), \tau(\bsy{d}) \in \cpf(G_H)$. 

    To show $\mathcal{O}' \in \mathcal{A}(G)_{H \to H^C}$, observe that for an arbitrary $\mathcal{O}_G \in \mathcal{A}(G)$, 
   \begin{eqnarray}
\label{eq:indegree_GA_sum}
    \sum_{ j \in H } \indeg_{\mathcal{O}_G}(j)& =& \overbrace{\sum_{\substack{ \{0,j\} \in E(G) \\ j \in V(H) }} wt_G(\{0,j\})}^{=C_1(\mathcal{O}_G)} + \overbrace{\sum_{\substack{ \{i,j\} \in E(G) \\ i,j \in V(H) }} wt_G(\{i,j\})}^{=C_2(\mathcal{O}_G)} + \overbrace{\sum_{\substack{ \{i,j\} \in E(G) \\ j \in V(H), \ i \in V(H^C) \\ \text{and } (i,j) \in \mathcal{O}_G }} wt_G(\{i,j\})}^{=C_3(\mathcal{O}_G)}  \nonumber \\
    &\geq &\sum_{ \{i,j\} \in E(G_H) } wt_G(\{i,j\}),  \label{Ineq:indeg} 
\end{eqnarray} 
as the last formula is the sum of $C_1(\mathcal{O}_G)$ and $C_2(\mathcal{O}_G)$. In particular, equation in \eqref{Ineq:indeg} holds if and only if $\mathcal{O}_G$ is in $\mathcal{A}(G)_{H \to H^C}$. 

Back to our problem. 
Since $\bsy{c}$ and $\sigma(\bsy{c})$ (resp. $\bsy{d}$ and $\tau(\bsy{d})$ ) have the same sum of their entries, 
$\sum_{j \in H} \indeg_{\mathcal{O}}(j) = \sum_{j \in H} \indeg_{\mathcal{O}'}(j)$. 
But $\mathcal{O} \in \mathcal{A}(G)_{H \to H^C}$. It follows that the equation  in \eqref{Ineq:indeg} holds for $\mathcal{O}$. Hence it  holds for $\mathcal{O}'$ and $\mathcal{O}' \in \mathcal{A}(G)_{H \to H^C}$.
\end{proof}


Before delving into the detailed discussion, we  present a list of key techniques frequently used in our case-by-case classification of graphs 
$G$ whose set of $G$-parking functions is invariant under the action of 
$\mathfrak{S}_p \times \mathfrak{S}_q$. 

\begin{enumerate}
    \item If $\cpf(G)$ is invariant under $\mathfrak{S}_p \times \mathfrak{S}_q$ for a graph $G$, all $A$-vertices (resp. $B$-vertices) must have the same set of obtainable indegrees under orientations in $\mathcal{A}(G)$. This follows from Lemma~\ref{lemma:invariance_maximal_GPF} and Proposition~\ref{burningalgorithm}. Consequently, if the values of the potential indegrees for $A$-vertices differ, or if one $A$-vertex has a larger set of potential indegrees than another, $\cpf(G)$ does not have the stated invariance property for the considered graph $G$.

    \item Since we are only concerned with orientations in the set $\mathcal{A}(G)$, a cut vertex $v$ of $G$ will always function as a unique source vertex to any (connected) subgraph of $G$ that contains $v$ but no vertices of the $0$-component of $G \setminus \{v\}$. Explicitly, if $C$ is any connected component of $G \setminus \{v\}$ not containing $0$, then for any $\mathcal{O} \in \mathcal{A}(G)$, $v$ is the unique source vertex in $G[\{v\} \cup V(C)]$ when $\mathcal{O}$ is restricted to that subgraph. This means that the possible indegrees of $v$ under orientations in $\mathcal{A}(G)$ are only dependent on the neighbors of $v$ in the $0$-component of $G \setminus \{v\}$.
    
    \item For the case that  $G_A$  is not connected, 
     we often use the following argument to show that the graph must have a simple symmetric structure: 
     Suppose the graph contains three vertices $u, v, w$ that are all in $A$ but are not in the same connected component of $G_A$.  
    Then we first construct an orientation $\mathcal{O} \in \mathcal{A}(G)$ such that $\indeg_\mathcal{O}(u)=\indeg_\mathcal{O}(v)=a$ and 
    $\indeg_\mathcal{O}(w)=b$, where $a\neq b$, then show that there does not exist an orientation $\mathcal{O}' \in \mathcal{A}(G)$ such that $\indeg_{\mathcal{O}'}(u)=\indeg_\mathcal{O}(w)=a$,  
    $\indeg_{\mathcal{O}'}(v)=b$,  and $\indeg_\mathcal{O}(x))= \indeg_{\mathcal{O}'}(x)$ for all vertices $x \neq u, v, w$.   
\end{enumerate}

\subsection{$G_A$ is connected}
\label{subsec:connect}

Unless otherwise specified, we will assume $G_A$ is connected in this entire subsection. All arguments symmetrically apply to the case where $G_B$ is connected. We begin with some notation and an important definition. We will say that a graph $G$ has type $\mathscr{G}_{_\mathrm{I}}(A, B)$ if $G$ satisfies the following:
\begin{enumerate}
    \item $V(G) = \{0\} \cup A \cup B$, where 0 is the unique root vertex of $G$, $A = \{1,2,\ldots,p\}$, and $B = \{p+1,p+2,\ldots,p+q\}$ for some $p,q \in \mathbb{Z}^+$;
    \item $G_A$ is connected;
    \item $\cpf(G)$ is invariant under the action of the product of symmetric groups $\mathfrak{S}_{p} \times \mathfrak{S}_{q}$.
\end{enumerate}

When $G_A$ is assumed to be connected, setting $H = G_A$ in Proposition~\ref{prop:subgraph_invariance} implies that $\cpf(G_A)$ is invariant under the action of $\mathfrak{S}_p$. 
Therefore we can get the structure of $G_A$ by Theorem~\ref{thm:gaydarov and hopkins}. 
We also need the following definition.

\begin{definition}
Given a graph $G$ and an equivalence relation $\sim$ on $V(G)$, the \textit{quotient graph} of a graph $G$ with respect to $\sim$, denoted by $G / \sim$, is the graph having vertex set equal to the quotient set $V(G) / \sim$ and edge set equal to $\{\{[i], [j]\} : \{i,j\} \in E(G)\}$, where $[i]$ indicates the equivalence class of vertex $i \in V(G)$. We also define a new weight function $wt_{\sim}: E(G / \sim) \to \mathbb{Z}^+$ to be associated with $G / \sim$ satisfying 
$$ wt_{\sim}(\{[i],[j]\}) = \sum_{u \in [i]} \sum_{v \in [j]} wt_G(\{u,v\}). $$
\end{definition}

An important ingredient in the proofs to follow will be the quotient graph of $G$ whose equivalence relation $\sim$ identifies all vertices in $A \cup \{0\}$. Using the notation $\bar{0}$ to represent the equivalence class of $A \cup \{0\}$, we will let $G / (A \cup \{0\})$ denote the graph $G / \sim$ with the added condition that $\bar{0}$ is the unique root vertex. Figure \ref{fig:quotient_graph} below depicts $G / (A \cup \{0\})$ for a certain graph $G$.


\begin{figure}[H]
\begin{center}
\begin{subfigure}[t]{.45\textwidth}
\begin{center}
\begin{tikzpicture}[line width=.7pt,scale=1]
\tikzstyle{vertex}=[draw,circle,fill=black,minimum size=3pt,inner sep=0pt]
\path (0,0) node[vertex, label=left:0] (v0) {};
\path (1,1) node[vertex, label=above:1] (v1) {};
\path (1,0) node[vertex, label=above left:2] (v2) {};
\path (1,-1) node[vertex, label=below:3] (v3) {};
\path (3,1) node[vertex,label=right:4] (v4) {};
\path (3,0) node[vertex,label=right:5] (v5) {};
\path (3,-1) node[vertex,label=right:6] (v6) {};
\path (4,0) node[vertex,label=right:7] (v7) {};

\draw[blue] (v0) edge (v1);
\draw[blue] (v0) edge (v2);
\draw[blue] (v0) edge (v3);
\draw[blue] (v1) edge (v2);
\draw[blue] (v1) edge [bend left=60] (v3);
\draw[blue] (v2) edge (v3);

\draw[red] (v1) edge (v4);
\draw[red] (v1) edge (v5);
\draw[red] (v2) edge (v5);
\draw[red] (v2) edge (v6);
\draw[red] (v3) edge (v5);
\draw[red] (v3) edge (v6);
\draw[red] (v5) edge (v6);
\draw[red] (v4) edge (v7);
\draw[red] (v6) edge (v7);

\end{tikzpicture}
\caption{A graph $G$ with $V(G) = \{0\} \cup A \cup B$, where $A = \{1,2,3\}$ and $B = \{4,5,6\}$.}
\label{fig:quotient_graph:sub1}
\end{center}
\end{subfigure}
\hspace{5mm}
\begin{subfigure}[t]{.45\textwidth}
\begin{center}
\begin{tikzpicture}[line width=.7pt,scale=1]
\tikzstyle{vertex}=[draw,circle,fill=black,minimum size=3pt,inner sep=0pt]
\path (0,0) node[vertex, label=left:$\bar{0}$] (v0) {};
\path (2,1) node[vertex,label=right:4] (v4) {};
\path (2,0) node[vertex,label=right:5] (v5) {};
\path (2,-1) node[vertex,label=right:6] (v6) {};
\path (3,0) node[vertex,label=right:7] (v7) {};

\draw[red] (v0) edge (v4);
\draw[red] (v0) edge (v5);
\draw[red] (v0) edge (v6);
\draw[red] (v5) edge (v6);
\draw[red] (v4) edge (v7);
\draw[red] (v6) edge (v7);

\end{tikzpicture}
\caption{The graph $G / (A \cup \{0\})$. For each $v \in B$, $wt_\sim(\{\bar{0}, v\}) = \sum_{u \in A} wt_G(\{u,v\})$.}
\label{fig:quotient_graph:sub2}
\end{center}
\end{subfigure}
\end{center}
\caption{An example of a graph $G$ and its quotient graph $G \setminus (A \cup \{0\})$.}
\label{fig:quotient_graph}
\end{figure}


\begin{proposition}
\label{prop:quotient_graph_surjection}
For any graph $G$ of type $\mathscr{G}_{_\mathrm{I}}(A, B)$, there exists a surjection from a subset of the maximal $G$-parking functions onto the set of maximal $(G / (A \cup \{0\}))$-parking functions.
\end{proposition}

\begin{proof} Let $\cMPF(G)_{A \to B}$ denote the subset of $\cMPF(G)$ consisting of all maximal $G$-parking functions which correspond to an orientation in $\mathcal{A}(G)_{A \to B}$ according to the bijection of Proposition \ref{burningalgorithm}. Define a map $f: \cMPF(G)_{A \to B} \to \cMPF(G / (A \cup \{0\}))$ by $(\bsy{a},\bsy{b}) \mapsto \bsy{b}$. Clearly $f$ is well-defined, according to the definition of $G / (A \cup \{0\})$. 
For any $\bsy{b} \in \cMPF(G / (A \cup \{0\})$, if $\mathcal{O}$ is the acyclic orientation of $G / (A \cup \{0\})$ with unique source $\bar{0}$ corresponding to $\bsy{b}$, we can extend $\mathcal{O}$ into an orientation $\mathcal{O}_G \in \mathcal{A}(G)_{A \to B}$ simply by fixing any valid acyclic orientation of $G_A$ with unique source 0.  Let  $\bsy{a}$ be the corresponding $G_A$-parking function. Then $(\bsy{a},\bsy{b})$ is an element of $\cMPF(G)_{A \to B}$. Hence, $f$ is the stated surjection. \end{proof}

\begin{corollary}
\label{cor:structure_quotient_graph}
If $G$ is a graph of type $\mathscr{G}_{_\mathrm{I}}(A, B)$, then $G / (A \cup \{0\})$ is a $c$-cycle or $c$-tree or $K_{1+q}^{c,d}$ for some $c,d \in \mathbb{Z}^+$.
\end{corollary}

\begin{proof} 
Let $(\bsy{a}, \bsy{b}) \in \cMPF_{A \to B}$ and $\tau \in \mathfrak{S}_q$. By the proof of \Cref{prop:subgraph_invariance}, the maximal $G$-parking function $(\bsy{a}, \tau(\bsy{a}))$ also corresponds to an orientation in $\mathcal{A}(G)_{A \to B}$. 
By \Cref{prop:quotient_graph_surjection}, we see that $\cMPF(G / (A \cup \{0\}))$ is invariant under the action of $\mathfrak{S}_q$, which implies by
\Cref{thm:gaydarov and hopkins}  that $G / (A \cup \{0\})$ is one of the three graphs in the corollary statement. 
\end{proof}

\subsubsection{$G$ has no cut vertices}
\label{subsubsection 6.2.1}

In the next three theorems, we use Corollary~\ref{cor:structure_quotient_graph} to classify all graphs of type $\mathscr{G}_{_\mathrm{I}}(A, B)$ which have no cut vertices.

\begin{theorem}
\label{thm:c-cycle and no CV}
Let $G$ be a graph of type $\mathscr{G}_{_\mathrm{I}}(A, B)$ such that $G$ has no cut vertices. Assume that either $G / (A \cup \{0\})$ is a $c$-cycle of length three with each $B$-vertex adjacent to a single (not necessarily the same) vertex in $A \cup \{0\}$, or $G / (A \cup \{0\})$ is a $c$-cycle of minimum length four, where $c \in \mathbb{Z}^+$. Then one of the following holds ($a,b \in \mathbb{Z}^+$):
\begin{enumerate}
    \item[(i)] $G$ is a $c$-cycle;
    \item[(ii)] $G$ is a cycle. $A=\{1\}$ and $wt_G(\{0,1\}) = a$ while $wt_G(\{i,j\}) = c$ for all other $\{i,j\} \in E(G)$;
    \item[(iii)] $A = \{1,2\}$ and $G$ is a cycle with possibly the chord $\{1, 2\}$. The neighbors of $0$ are exactly $1,2$ with 
     $wt_G(\{0,1\}) = wt_G(\{0,2\}) = a$, $wt_G\{1,2\}=b$ if $\{1,2\} \in E(G)$. For all other $\{i,j\} \in E(G)$, $wt_G(\{i,j\}) = c$. (cf. Figure~\ref{fig:c-cycle_examples:sub1})
\end{enumerate}
\end{theorem}

\noindent 
\textsc{Remark.} 
We include the case where $G / (A \cup \{0\})$ is a length-three $c$-cycle with each $B$-vertex adjacent to only one vertex in $A \cup \{0\}$ in the hypotheses of Theorem~\ref{thm:c-cycle and no CV} since graphs $G$ satisfying this condition are of a similar structure to graphs for which $G / (A \cup \{0\})$ is a $c$-cycle of minimum length four. The case where $G / (A \cup \{0\})$ is a length-three $c$-cycle with at least one $B$-vertex adjacent to multiple vertices in $A \cup \{0\}$ will be considered in Theorem~\ref{thm:complete and no CV}, when we assume $G / (A \cup \{0\})$ is $K^{c,d}_{q+1}$.

\medskip

\noindent \textit{Proof.} By hypothesis, $|B| \geq 2$. Without loss of generality, assume $p+1$ and $p+q$ are the two vertices of $B$ adjacent to the root $\bar{0}$ in $G / (A \cup \{0\})$, with $p+k$ adjacent to $p+k+1$ for each $1 \leq k \leq q-1$. Observe that we cannot have both $p+1$ and $p+q$ adjacent to the same vertex in $A \cup \{0\}$ and no other vertices in $A \cup \{0\}$, else $G$ would have a cut vertex. Thus, the only two possibilities under the theorem's hypotheses are that $p+1$ and $p+q$ are each adjacent to exactly one vertex in $A \cup \{0\}$ but not the same vertex, or $G / (A \cup \{0\})$ has minimum length four and at least one of $p+1$ and $p+q$ is adjacent to more than one vertex in $A \cup \{0\}$. 

\vspace{3mm}
\fbox{\textbf{Case 1.}} \textit{The vertices $p+1$ and $p+q$ are each adjacent to exactly one vertex in $A \cup \{0\}$ but not the same vertex.} Assume $p+1$ is adjacent to $i_1 \in A \cup \{0\}$ and $p+q$ is adjacent to $i_q \in A \cup \{0\}$ ($i_1 \neq i_q$). We know $A \subset \mathcal{I}(G)$ since $G$ has no cut vertices, and so by Lemma~\ref{lemma:non_cut_vertex_equal_degrees} all vertices in $A$ have the same degree in $G$.

If $|A|=1$ so that $A = \{1\}$, then clearly we can assume $i_1 = 1$ and $i_q = 0$, and $G$ must be a cycle with $wt_G(\{0,1\}) = a$ and all other edges of weight $c$. It $|A| = 2$ so that $A = \{1,2\}$, $G_A$ must be either a two-leaf $a$-star or the triangle $K_3^{a,b}$ (cf. Figure \ref{fig:c-cycle_examples:sub1}). Since both vertices of $A$ must have the same degree in $G$, we can assume $i_1 = 1$ and $i_q = 2$. All edges of $G$ not in $E(G_A)$ will have weight $c$.

If $|A| \geq 3$ and $G_A$ is an $a$-cycle or the complete graph $K_{p+1}^{a,b}$, then vertices $i_1$ and $i_q$ will not have the same degree as other $A$-vertices in $G$, violating Lemma~\ref{lemma:non_cut_vertex_equal_degrees}. Thus, $G_A$ must be an $a$-tree if $|A| \geq 3$. In fact, it is an easy consequence of Lemma~\ref{lemma:non_cut_vertex_equal_degrees} that $G_A$ must be a $c$-path with $i_1$ and $i_q$ its two leaf vertices, such that $G$ is a $c$-cycle.

\vspace{3mm}
\fbox{\textbf{Case 2.}} \textit{$G / (A \cup \{0\})$ has minimum length four, and at least one of $p+1$ and $p+q$ is adjacent to more than one vertex in $A \cup \{0\}$. } Assume $p+1$ is adjacent to more than one vertex in $A \cup \{0\}$. Since $G / (A \cup \{0\})$ is a $c$-cycle, we know $\sum_{i \in A \cup \{0\}} wt_G(\{i,p+1\}) = c$, meaning every edge of $E(G)$ between $p+1$ and a vertex in $A \cup \{0\}$ has weight strictly less than $c$. Note that $p+1$ must be adjacent to at least one vertex in $A$. We will show that there must exist an orientation in $\mathcal{A}(G)$ under which $p+1$ has indegree less than $c$.

If $G_A$ is either an $a$-cycle or the complete graph $K_{p+1}^{a,b}$, then we can fix a vertex $i \in A$ which is adjacent to $p+1$ and find an orientation $\mathcal{O}_A \in \mathcal{A}(G_A)$ under which vertex $i$ is a sink. Otherwise, if $G_A$ is an $a$-tree, we can choose the vertex $i \in A$ to have maximal distance from the root $0$ in $G$ among all vertices in the set $\{u \in A : \{u, p+1\} \in E(G)\}$; then let $\mathcal{O}_A$ be the unique orientation in $\mathcal{A}(G_A)$. Now extend $\mathcal{O}_A$ into an orientation $\mathcal{O}$ of $G$ by doing the following: direct all edges incident to $p+q$ towards $p+q$, direct every edge $\{p+k, p+k+1\}$ ($1 \leq k \leq q-1$) towards $p+k+1$, and direct all edges between $p+1$ and a vertex in $A \cup \{0\}$ towards $p+1$ with the exception of the edge $\{i,p+1\}$, which should be directed towards $i$. See Figure \ref{fig:c-cycle_examples:sub2} for a representation of $\mathcal{O}$. Then $\mathcal{O} \in \mathcal{A}(G)$ and $\indeg_{\mathcal{O}}(p+1) < c$. However, since $|B| \geq 3$, there is at least one vertex $p+k$, $2 \leq k \leq q-1$, which can never have indegree less than $c$ under any orientation in $\mathcal{A}(G)$ since 0 must be the unique source. Therefore, Case 2 is not possible. \hfill \qedsymbol \\


\begin{figure}[H]
\begin{center}
\begin{subfigure}[t]{.45\textwidth}
\begin{center}
\begin{tikzpicture}[line width=.7pt,scale=1]
\tikzstyle{vertex}=[draw,circle,fill=black,minimum size=3pt,inner sep=0pt]
\path (0,0) node[vertex, label=above:0] (v0) {};
\path (-1,-1) node[vertex, label=left:1] (v1) {};
\path (1,-1) node[vertex, label=right:2] (v2) {};
\path (-2,-2) node[vertex,label=below:3] (v3) {};
\path (-1,-2) node[vertex,label=below:4] (v4) {};
\path (0,-2) node (ell) {$\cdots$};
\path (1,-2) node[vertex,label=below:$q+1$] (vq+1) {};
\path (2,-2) node[vertex,label=below:$q+2$] (vq+2) {};

\draw[blue] (v0) -- (v1) node [midway, above=0.5mm] {$a$};
\draw[blue] (v0) edge (v2);

\draw[Green] (v1) -- (v2) node [midway,above=0.01mm] {$b$};

\draw[red] (v1) -- (v3) node [midway, left=0.4mm] {$c$};
\draw[red] (v3) edge (v4);
\draw[red] (v4) edge (ell);
\draw[red] (ell) edge (vq+1);
\draw[red] (vq+1) edge (vq+2);
\draw[red] (vq+2) edge (v2);

\end{tikzpicture}
\caption{A graph $G$ described in Theorem~\ref{thm:c-cycle and no CV}(iii).}
\label{fig:c-cycle_examples:sub1}
\end{center}
\vspace{4mm}
\end{subfigure}
\hspace{5mm}
\begin{subfigure}[t]{.45\textwidth}
\begin{center}
\begin{tikzpicture}[line width=.7pt,scale=1]
\tikzstyle{vertex}=[draw,circle,fill=black,minimum size=3pt,inner sep=0pt]
\path (0,0) node[vertex, label=above:0] (v0) {};
\path (0,-0.75) node (ell_A) {\rvdots};
\path (0,-1.5) node[vertex, label=right:$i$] (vi) {};
\path (-1.25,-.75) node (v1) {};
\path (1.25,-.75) node (v2) {};
\path (0,-.6) node[ellipse,minimum width = 40mm,minimum height = 25mm,draw,dashed,blue,label = above right:$G_A$] (G_A) {};
\path (-3,-3) node[vertex, label=below:$p+1$] (vp+1) {};
\path (-1.5,-3) node[vertex, label=below:$p+2$] (vp+2) {};
\path (0,-3) node (ell_B) {$\cdots$};
\path (1.5,-3) node[vertex, label=below:$p+q-1$] (vp+q-1) {};
\path (3,-3) node[vertex, label=below:$p+q$] (vp+q) {};

\draw[-Stealth,blue] (v0) edge (ell_A);
\draw[-Stealth,blue] (ell_A) edge (vi);

\draw[-Stealth,red] (vp+1) edge (vi);
\draw[-{Stealth[].Stealth[]Stealth[]},red] (v1) edge (vp+1);
\draw[-Stealth,red] (vp+1) edge (vp+2);
\draw[-Stealth,red] (vp+2) edge (ell_B);
\draw[-Stealth,red] (ell_B) edge (vp+q-1);
\draw[-Stealth,red] (vp+q-1) edge (vp+q);
\draw[-{Stealth[].Stealth[]Stealth[]},red] (v2) edge (vp+q);
 
\end{tikzpicture}
\caption{An $\mathcal{O} \in \mathcal{A}(G)$ in which $\indeg_{\mathcal{O}}(p+1)$ has smaller indegree than is possible for the vertices in $B \setminus \{p+1,p+q\}$. A triple headed arrow indicates one or more edges in that direction.}
\label{fig:c-cycle_examples:sub2}
\end{center}
\end{subfigure}
\end{center}
\caption{A graph (on left) and an orientation (on right) related to Theorem~\ref{thm:c-cycle and no CV}.}
\label{fig:c-cycle_examples}
\end{figure}


\begin{theorem}
\label{thm:c-tree and no CV}
Let $G$ be a graph of type $\mathscr{G}_{_\mathrm{I}}(A, B)$. If $G$ has no cut vertices and $G / (A \cup \{0\})$ is a $c$-tree for some $c \in \mathbb{Z}^+$, then $G$ is characterized as follows:
\begin{itemize}
    \item $G_A$ is either an $a$-star with root 0 or equal to $K_{p+1}^{a,b}$ for $a,b \in \mathbb{Z}^+$; and
    \item every vertex of $B$ is adjacent to every vertex of $A$ via a edge of a constant weight $c \in \mathbb{Z}^+$; and
    \item either $G_B$ is a $d$-star with root 0 for $d \in \mathbb{Z}^+$ or $G_B$ consists of all singleton vertices.
\end{itemize}
\end{theorem}

\begin{proof}  We first observe that if $G / (A \cup \{0\})$ is a $c$-tree and $G$ has no cut vertices, every vertex $j \in B \subset \mathcal{I}(G)$ must be adjacent to the root $\bar{0}$ in $G / (A \cup \{0\})$. Thus, $G / (A \cup \{0\})$ must be a $c$-star with root $\bar{0}$. We wish to show the following: (1) every vertex of $A$ must be adjacent to every vertex of $B$; (2) $G_A$ is either an $a$-star or $K_{p+1}^{a,b}$; (3) every edge between a vertex of $A$ and a vertex of $B$ has the same positive weight; 
and (4) if one vertex of $B$ is adjacent to the root 0, then all vertices of $B$ are adjacent to the root 0 via edges of equal weight.

\vspace{3mm}
\fbox{\textbf{Claim 1.}} \textit{Every vertex of $A$ is adjacent to every vertex of $B$.}

\textit{Proof of Claim 1.} Suppose there exist vertices $i \in A$ and $j \in B$ such that $\{i,j\} \notin E(G)$. Then we can find an orientation $\mathcal{O} \in \mathcal{A}(G)$ with $\indeg_{\mathcal{O}}(i) = deg_G(i)$ and $\indeg_{\mathcal{O}}(j) = \deg_G(j)$ (fix an orientation $\mathcal{O}_A \in \mathcal{A}(G_A)$ for which $\indeg_{\mathcal{O}_A}(i) = \deg_{G_A}(i)$ and extend it into an orientation $\mathcal{O} \in \mathcal{A}(G)$ by directing all edges incident to $i$ and a vertex of $B$ towards $i$, and then directing all remaining edges towards the $B$-vertex). Let $i' \in A$ be adjacent to $j$. Then by Lemma~\ref{lemma:invariance_maximal_GPF}, there should exist an $\mathcal{O}' \in \mathcal{A}(G)$ with $\indeg_{\mathcal{O}'}(i')  =\deg_G(i)=\deg_G(i')$ and $\indeg_{\mathcal{O}'}(j) = deg_G(j)$. But this is clearly not possible. 
We conclude that every vertex of $A$ is adjacent to every vertex of $B$. $\blacksquare$

\vspace{3mm}
\fbox{\textbf{Claim 2.}} \textit{$G_A$ is either an $a$-star or $K_{p+1}^{a,b}$.}

\textit{Proof of Claim 2.} If $0$ is not adjacent to some vertex of $A$, then $G_A$ is either an $a$-tree but not $a$-star, or an $a$-cycle with length at least 4. Let us consider these two cases separately.

\indent \indent \fbox{\textit{Case 2a.}} \textit{$G_A$ is an $a$-tree but not $a$-star.} Then there must be a leaf $\ell \in A$ with $\ell \notin N_G(0)$. Thus, there are vertices $i_1, i_2, \ldots, i_n \in A$ $(n \geq 1)$ such that $(0, i_1, i_2, \ldots, i_n, \ell)$ forms a directed path from 0 to $\ell$. Consider the orientation $\mathcal{O} \in \mathcal{A}(G)$ satisfying the following: all edges in the tree $G_A$ are directed towards the leaf vertices, all edges incident to a vertex in the set $A \setminus \{i_1, i_2, \ldots, i_n\}$ and a $B$-vertex are directed towards the $A$-vertex, and all edges incident to a vertex in the set $\{0,i_1, i_2, \ldots, i_n\}$ and a $B$-vertex are directed towards the $B$-vertex (cf. Figure \ref{fig:c-tree_example_1}). 


\begin{figure}[H]
\begin{center}
\begin{tikzpicture}[line width=.7pt,scale=1]
\tikzstyle{vertex}=[draw,circle,fill=black,minimum size=3pt,inner sep=0pt]
\path (0,1.25) node[vertex, label=left:0] (v0) {};
\path (-0.5,0.25) node[vertex] (v1) {};
\path (1, 1.25) node[vertex, label=above:$i_1$] (vi1) {};
\path (2,1.25 ) node (ell_A) {$\cdots$};
\path (3,1.25) node[vertex, label=above:$i_n$] (vin) {};
\path (4,1.25) node[vertex, label=above:$\ell$] (l) {};
\path (1.5,-1.25) node[vertex, label=below:$p+1$] (vp+1) {};
\path (2.5,-1.25) node (ell_B) {$\cdots$};
\path (3.5,-1.25) node[vertex, label=below:$p+q$] (vp+q) {};

\draw[-Stealth, blue] (v0) edge (vi1);
\draw[-Stealth, blue] (vi1) edge (ell_A);
\draw[-Stealth, blue] (ell_A) edge (vin);
\draw[-Stealth, blue] (vin) edge (l);
\draw[-Stealth, blue] (v0) edge (v1);

\draw[-Stealth, red] (vi1) edge (vp+1);
\draw[-Stealth, red] (vi1) edge (vp+q);
\draw[-Stealth, red] (vin) edge (vp+1);
\draw[-Stealth, red] (vin) edge (vp+q);
\draw[-Stealth, red] (vp+1) edge (l);
\draw[-Stealth, red] (vp+q) edge (l);
\draw[-Stealth, red] (vp+1) edge[bend left=15]  (v1);
\draw[-Stealth, red] (vp+q) edge[bend left=15] (v1);

\end{tikzpicture}
\caption{The orientation $\mathcal{O}$ in Case 2a. Note that if 0 is not adjacent to any vertices of $B$, then 0 must be adjacent to at least two vertices of $A$, (otherwise $i_1$ is a cut vertex .}
\label{fig:c-tree_example_1}
\end{center}
\end{figure}


Observe that the indegree of every leaf vertex of $A$ with respect to $\mathcal{O}$ is maximal; in particular, $\indeg_{\mathcal{O}}(\ell) = deg_G(\ell)$. By Lemma~\ref{lemma:invariance_maximal_GPF}, there should be another orientation $\mathcal{O}' \in \mathcal{A}(G)$ for which $\indeg_{\mathcal{O}'}(\ell) = \indeg_{\mathcal{O}}(i_n) = a$, $\indeg_{\mathcal{O}'}(i_n) = \indeg_{\mathcal{O}}(\ell) = \deg_G(\ell) = \deg_G(i_n)$, and $\indeg_{\mathcal{O}'}(u) = \indeg_{\mathcal{O}}(u)$ for all other $u \in V(G) \setminus \{i_n, \ell\}$. The last condition implies that any edge incident to a  vertex in $A-\{i_1, \dots, i_n, \ell\}$ and a $B$-vertex must be directed towards the $A$-vertex in $\mathcal{O}'$. 
Since $\ell$ has zero indegree when $\mathcal{O}'$ is restricted to $G_A$, there must be at least one $j \in B$ with $(j,\ell) \in \mathcal{O}'$. 
It follows that 

$$ \indeg_{\mathcal{O}}(j) = \sum_{i \in \{0,i_1,\ldots,i_n\}} wt_G(\{i,j\}) \quad\quad \text{and}  \quad\quad \indeg_{\mathcal{O}'}(j) = \sum_{i \in \{0,i_1,\ldots,i_{n-1}\}} wt_G(\{i,j\}),$$
meaning $\indeg_{\mathcal{O}}(j) \neq \indeg_{\mathcal{O}'}(j)$, which contradicts the definition of $\mathcal{O}'$.

\indent \indent \fbox{\textit{Case 2b.}} \textit{$G_A$ is an $a$-cycle of length at least 4.}\  In this case, $p \geq 3$. Assume the structure of $G_A$ is such that 0 is adjacent to vertices 1 and $p$ and vertex $i$ is adjacent to vertex $i+1$ for $1 \leq i < p$, let $\mathcal{O} \in \mathcal{A}(G)$ satisfy the following: $G_A$ contains the directed path $(1, 2, \ldots, p)$, all edges between 1 and a $B$-vertex are directed towards the $B$-vertex, and all edges between a vertex in $\{2, \ldots, p\} \subset A$ and a $B$-vertex are directed towards the $A$-vertex (cf. Figure \ref{fig:c-tree_example_2}. Edges between $0$ and $B$ are omitted since they are directed towards $B$ in any orientation of $\mathcal{A}(G)$. ). 


\begin{figure}[H]
\begin{center}
\begin{tikzpicture}[line width=.7pt,scale=0.95]
\tikzstyle{vertex}=[draw,circle,fill=black,minimum size=3pt,inner sep=0pt]
\path (-1,0) node[vertex, label=left:0,] (v0) {};
\path (1,1.5) node[vertex, label=above:1] (v1) {};
\path (1,0.75) node[vertex, label=left:2] (v2) {};
\path (1,0) node (ell_A) {\rvdots};
\path (1,-0.75) node[vertex, label=left:$p-1$] (vp-1) {};
\path (1,-1.5) node[vertex, label=below:$p$] (vp) {};
\path (3.5,0.75) node[vertex, label=right:$p+1$] (vp+1) {};
\path (3.5,0) node (ell_B) {\rvdots};
\path (3.5,-0.75) node[vertex, label=right:$p+q$] (vp+q) {};

\draw[-Stealth, blue] (v0) edge[bend left=15] (v1);
\draw[-Stealth, blue] (v1) edge (v2);
\draw[-Stealth, blue] (v2) edge (ell_A);
\draw[-Stealth, blue] (ell_A) edge (vp-1);
\draw[-Stealth, blue] (vp-1) edge (vp);
\draw[-Stealth, blue] (v0) edge[bend right=15] (vp);

\draw[-Stealth, red] (v1) edge[bend left=15] (vp+1);
\draw[-Stealth, red] (v1) edge[bend left=15] (vp+q);
\draw[-Stealth, red] (vp+1) edge (v2);
\draw[-Stealth, red] (vp+1) edge (vp-1);
\draw[-Stealth, red] (vp+1) edge (vp);
\draw[-Stealth, red] (vp+q) edge[bend left=15] (v2);
\draw[-Stealth, red] (vp+q) edge[bend left=15] (vp-1);
\draw[-Stealth, red] (vp+q) edge[bend left=15] (vp);

\end{tikzpicture}
\caption{The orientation $\mathcal{O}$ in Case 2b.}
\label{fig:c-tree_example_2}
\end{center}
\end{figure}


Note that $\indeg_{\mathcal{O}}(1) =a$, $\indeg_{\mathcal{O}}(p) = deg_G(p)$, and $\indeg_{\mathcal{O}}(i) = deg_G(i)-a$ for all $i \in A \setminus \{1,p\}$ (with $deg_G(i)$ constant for all $A$-vertices). Then by Lemma~\ref{lemma:invariance_maximal_GPF}, there must be an $\mathcal{O}' \in \mathcal{A}(G)$ such that $\indeg_{\mathcal{O}'}(1) = \indeg_{\mathcal{O}}(2) = deg_G(1)-a$, $\indeg_{\mathcal{O}'}(2) = a$, and $\indeg_{\mathcal{O}'}(i) = \indeg_{\mathcal{O}}(i)$ for all $i \in A \setminus \{1,2\}$. However, to maintain acyclicity, it must be the case in $\mathcal{O}'$ that the edge $\{1,2\}$ is directed towards 1  and at least one edge $\{2,j\}$ for some $j \in B$ is directed towards 2. Consequently, the edge $\{1,j\}$ must be directed towards 1 to maintain acyclicity, which then makes $j$ a source vertex in $\mathcal{O}'$, a contradiction.

Therefore, the root 0 must be adjacent to all vertices of $A$, and $G_A$ is either an $a$-star or equal to $K_{p+1}^{a,b}$ for some $a,b \in \mathbb{Z}^+$. $\blacksquare$

\vspace{3mm}
\fbox{\textbf{Claim 3.}} \textit{All edges between $A$-vertices and $B$-vertices have the same weight.}

\textit{Proof of Claim 3.} Based on the possible structures of $G_A$ given by Claim 2, there must be a constant $\alpha$ such that $deg_{G_A}(i) = \alpha$ for all $i \in A$. Moreover, $\alpha = \max_{\mathcal{O}_A \in \mathcal{A}(G_A)} \indeg_{\mathcal{O}_A}(i)$ for all $i \in A$. We also know that $deg_G(i)$ is constant for all $i \in A$ since $A \subset \mathcal{I}(G)$, and thus the sum $N = \sum_{j \in B} wt_G(\{i,j\})$ must also be constant for all $i \in A$. Hence, to show that all edges between $A$-vertices and $B$-vertices have the same weight, it is sufficient to fix an arbitrary vertex $v \in B$ and show that all edges incident to $v$ and an $A$-vertex have the same weight. 

Let $\mathcal{O} \in \mathcal{A}(G)$ satisfy the following: $(v,1) \in \mathcal{O}$, $(i,v) \in \mathcal{O}$ for all $i \in A \setminus \{1\}$, $(i,j) \in \mathcal{O}$ for all $i \in A$ and $j \in B \setminus \{v\}$, and $\indeg_{\mathcal{O}}(1) = \alpha + wt_G(\{1,v\})$. Then by Lemma~\ref{lemma:invariance_maximal_GPF}, for each $i \in A$, there should exist an orientation $\mathcal{O}_i \in \mathcal{A}(G)$ for which $\indeg_{\mathcal{O}_i}(i) = \indeg_{\mathcal{O}}(1) = \alpha + wt_G(\{1,v\})$ and $\indeg_{\mathcal{O}_i}(j) = \indeg_{\mathcal{O}}(j)$ for all $j \in B$. Specifically, $\indeg_{\mathcal{O}_i}(v) = \indeg_{\mathcal{O}}(v) = deg_G(v) - wt_G(\{1,v\})$ and $\indeg_{\mathcal{O}_i}(j) = \indeg_{\mathcal{O}}(j) = deg_G(j)$ for all $j \in B \setminus \{v\}$. This implies that $\indeg_{\mathcal{O}_i}(i) = \alpha + wt_G(\{i,v\})$, and hence $wt_G(\{i,v\}) = wt_G(\{1,v\})$ for all $i \in A$. $\blacksquare$

\vspace{3mm}
\fbox{\textbf{Claim 4.}} \textit{If $N_G(0) \cap B \neq \emptyset$, then all $B$-vertices are adjacent to $0$ via edges of equal weight.}

\textit{Proof of Claim 4.} This statement is an immediate consequence of the fact that $deg_G(j)$ is constant for all $j \in B$ and Claim 3. $\blacksquare$

\vspace{3mm}
The theorem follows by Claims 1--4. 
\end{proof} 

\medskip

\begin{theorem}
\label{thm:complete and no CV}
Let $G$ be a graph of type $\mathscr{G}_{_\mathrm{I}}(A, B)$ such that $G$ has no cut vertices. Assume that either $G / (A \cup \{0\})$ is the triangle $K_{3}^{e,d}$ with at least one $B$-vertex adjacent to more than one vertex in $A \cup \{0\}$, or $G / (A \cup \{0\})$ is equal to $K_{q+1}^{e,d}$ for $q \geq 3$, where $d,e \in \mathbb{Z}^+$. Then $G$ is characterized as follows:
\begin{itemize}
    \item $G_A$ is either an $a$-star with root 0 or equal to $K_{p+1}^{a,b}$ for $a,b \in \mathbb{Z}^+$; and
    \item every vertex of $B$ is adjacent to every vertex of $A$ via an equally weighted edge; and
    \item either $G_B$ is equal to $K_{q+1}^{c,d}$ for $c \in \mathbb{Z}^+$, or 
    $G[B]$ is equal to $K_q^d$ and there is no edge between  0 and $B$-vertices.  
\end{itemize}
\end{theorem}

\begin{proof} We wish to show that the same four statements proved in Theorem~\ref{thm:c-tree and no CV} hold. 

\vspace{3mm}
\fbox{\textbf{Claim 1.}} \textit{Every vertex of $A$ is adjacent to every vertex of $B$.}

\textit{Proof of Claim 1.} The argument used to prove Claim 1 in Theorem~\ref{thm:c-tree and no CV} applies, with one additional condition: in the orientation $\mathcal{O}$, edges in $G[B] = K_q^d$ incident to vertex $j$ will also be directed towards $j$ to ensure $\indeg_{\mathcal{O}}(j) = deg_G(j)$. $\blacksquare$

\vspace{3mm}
\fbox{\textbf{Claim 2.}} \textit{$G_A$ is either an $a$-star or $K_{p+1}^{a,b}$.}

\textit{Proof of Claim 2.} The argument used to prove Claim 2 in Theorem~\ref{thm:c-tree and no CV} applies in the case where $G_A$ is an $a$-tree but not $a$-star, with the added condition that the orientation $\mathcal{O}$ should also be acyclic when restricted to the graph $G_B$. 

A new argument is needed in the case where $G_A$ is assumed to be an $a$-cycle with $p \geq 3$. Assume the structure of $G_A$ is such that 0 is adjacent to vertices 1 and $p$, and vertex $i$ is adjacent to vertex $i+1$ for $1 \leq i < p$. Let $N = \sum_{j \in B} wt_G(\{i,j\})$, which is constant for all $i \in A$ since $A \subset \mathcal{I}(G)$. First, observe that if there exists at least one vertex of $B$ which is incident to 0, then the set $\mathcal{A}(G)_{B \to A} \subset \mathcal{A}(G)$ is non-empty. 
Construct  $\mathcal{O} \in \mathcal{A}(G)_{B \to A}$ as follows: 
\begin{enumerate}
    \item Fix an order on the elements in $B$  that starts with the neighbors of $0$. For $i,j \in B$, the edge $\{i,j\}$ is directed towards $j$ whenever $i$ is ordered  before $i$.
    \item $(2,1) \in \mathcal{O}$ and $(2,3) \in \mathcal{O}$.
\end{enumerate}
Then $\indeg_{\mathcal{O}}(2) = N$ and $\indeg_{\mathcal{O}}(1) = 2a + N$. Lemma~\ref{lemma:invariance_maximal_GPF} implies the existence of some $\mathcal{O'} \in \mathcal{A}(G)$ such that $\indeg_{\mathcal{O'}}(1) = N$, $\indeg_{\mathcal{O'}}(2) = 2a + N$, and $\indeg_{\mathcal{O'}}(i) = \indeg_{\mathcal{O}}(i)$ for all $3 \leq i \leq p+q$. But since such an $\mathcal{O'}$ would necessarily be contained in $\mathcal{A}(G)_{B \to A}$, we would have $\indeg_{\mathcal{O'}}(1) \geq a + N$. Therefore, $G_A$ cannot be an $a$-cycle if at least one vertex in $B$ is adjacent to the root 0.

Otherwise, assume no vertices of $B$ are adjacent to 0. Choose an orientation $\mathcal{O} \in \mathcal{A}(G)$ satisfying the following: $G_A$ contains the directed path $(1, 2, \ldots, p)$; $(j,p-1),(j,p) \in \mathcal{O}$ for $j \in B$; all edges incident to a vertex in $A \setminus \{p-1,p\}$ and a $B$-vertex are directed towards the $B$-vertex, and the edges of $G[B]$ are directed in any valid acyclic way (cf. Figure \ref{fig:complete_example_1}).


\begin{figure}[H]
\begin{center}
\begin{tikzpicture}[line width=.7pt,scale=0.95]
\tikzstyle{vertex}=[draw,circle,fill=black,minimum size=3pt,inner sep=0pt]
\path (-1,0) node[vertex, label=left:0,] (v0) {};
\path (1,1.5) node[vertex, label=above:1] (v1) {};
\path (1,0.75) node[vertex, label=left:2] (v2) {};
\path (1,0) node (ell_A) {\rvdots};
\path (1,-0.75) node[vertex, label=left:$p-1$] (vp-1) {};
\path (1,-1.5) node[vertex, label=below:$p$] (vp) {};
\path (3.5,0.75) node[vertex, label=above right:$p+1$] (vp+1) {};
\path (3.5,0) node (ell_B) {\rvdots};
\path (3.5,-0.75) node[vertex, label=below right:$p+q$] (vp+q) {};

\draw[-Stealth, blue] (v0) edge[bend left=15] (v1);
\draw[-Stealth, blue] (v1) edge (v2);
\draw[-Stealth, blue] (v2) edge (ell_A);
\draw[-Stealth, blue] (ell_A) edge (vp-1);
\draw[-Stealth, blue] (vp-1) edge (vp);
\draw[-Stealth, blue] (v0) edge[bend right=15] (vp);

\draw[-Stealth, red] (v1) edge[bend left=15] (vp+1);
\draw[-Stealth, red] (v1) edge[bend left=15] (vp+q);
\draw[-Stealth, red] (v2) edge (vp+1);
\draw[-Stealth, red] (vp+1) edge (vp-1);
\draw[-Stealth, red] (vp+1) edge (vp);
\draw[-Stealth, red] (v2) edge[bend right=15] (vp+q);
\draw[-Stealth, red] (vp+q) edge[bend left=15] (vp-1);
\draw[-Stealth, red] (vp+q) edge[bend left=15] (vp);
\draw[Green] (vp+1) edge[bend left=30] (vp+q);

\end{tikzpicture}
\caption{The orientation $\mathcal{O}$.}
\label{fig:complete_example_1}
\end{center}
\end{figure}


By our definition, the ordered indegree sequence for the $A$-vertices under $\mathcal{O}$ is $(a, a, \ldots, a, a + N, 2a + N)$. By Lemma~\ref{lemma:invariance_maximal_GPF}, there should exist some $\mathcal{O}' \in \mathcal{A}(G)$ whose ordered indegree sequence for the $A$-vertices is $(a + N, a, a, \ldots, a, 2a + N)$. Then under $\mathcal{O}'$, all edges of $G$ incident to $p$ are directed towards $p$, and either $G_A$ still contains the directed path $(1, 2, \ldots, p)$ or at least one vertex in $A \setminus \{p\}$ has zero indegree under $\mathcal{O}'$ restricted to $G_A$. We will consider these two cases separately.

\indent \indent \fbox{\textit{Case 3a.}} \textit{The directed path $(1,2, \ldots, p)$ is in $\mathcal{O}'$.} This means all edges incident to a vertex in $A \setminus \{1,p\}$ and a $B$-vertex are directed towards the $B$-vertex in $\mathcal{O}'$. Also, all edges $\{1,j\}$ for $j \in B$ are directed towards 1. Clearly $\mathcal{O}'$ has a cycle in this case, a contradiction.

\indent \indent \fbox{\textit{Case 3b.}} \textit{At least one vertex in $A \setminus \{1,p\}$ has zero indegree when $\mathcal{O}'$ is restricted to $G_A$.} Suppose there is exactly one such vertex, $i \in A \setminus \{1,p\}$. Since $\indeg_{\mathcal{O}'}(i) = a$, at least one edge incident to $i$ and a $B$-vertex is directed towards $i$. Furthermore, since all vertices in $A \setminus \{1,p\}$ have indegree $a$ under $\mathcal{O}'$, $\mathcal{O}'$ must contain the directed paths $(i, i-1, \ldots, 2, 1)$ and $(i, i+1, \ldots, p-1, p)$, and all edges incident to a vertex in the set $\{2, 3, \ldots, i-1, i+1, \ldots, p-1\}$ and a $B$-vertex are directed towards the $B$-vertex. 

If $p > 3$, then $i \in \{2, 3, \ldots, p-1\}$. Since $\indeg_{\mathcal{O}'}(i) = a$, there is at least one vertex $j \in B$ such that $(j,i) \in \mathcal{O}'$. But then there must be some vertex $i' \in \{2, \ldots, i-1, i+1, p-1\}$ for which $\mathcal{O}'$ contains the cycle $(i,i',j,i)$, a contradiction. Therefore, $p=3$ and $i=2$.


\begin{figure}[H]
\begin{center}
\begin{tikzpicture}[line width=.7pt,scale=1]
\tikzstyle{vertex}=[draw,circle,fill=black,minimum size=3pt,inner sep=0pt]
\path (0,0) node[vertex, label=above:0] (v0) {};
\path (-3,-1) node[vertex, label=left:1] (v1) {};
\path (-2,-1) node (ell_A1) {$\cdots$};
\path (-1,-1) node[vertex, label=above:$i'$] (vi') {};
\path (0,-1) node (ell_A2) {$\cdots$};
\path (1,-1) node[vertex, label=above:$i$] (vi) {};
\path (2,-1) node (ell_A3) {$\cdots$};
\path (3,-1) node[vertex, label=right:$p$] (vp) {};
\path (0,-2.5) node[vertex, label=below:$j$] (vj) {};

\draw[-Stealth, blue] (v0) edge[bend right=15] (v1);
\draw[Stealth-, blue] (v1) edge (ell_A1);
\draw[Stealth-, blue] (ell_A1) edge (vi');
\draw[Stealth-, blue] (vi') edge (ell_A2);
\draw[Stealth-, blue] (ell_A2) edge (vi);
\draw[-Stealth, blue] (vi) edge (ell_A3);
\draw[-Stealth, blue] (ell_A3) edge (vp);
\draw[-Stealth, blue] (v0) edge[bend left=15] (vp);

\draw[-Stealth, red] (vj) edge (vi);
\draw[-Stealth, red] (vi') edge (vj);

\path (0,-2.5) node[ellipse, minimum width = 40mm,minimum height = 15mm,draw,dashed,Green,label = below right:$B$] (B) {};

\end{tikzpicture}
\caption{The orientation $\mathcal{O}'$ in Case 3b.}
\label{fig:complete_example_2}
\end{center}
\end{figure}


Let $B' = \{j \in B : (j,2) \in \mathcal{O}'\}$; then $\sum_{j \in B'} wt_G(\{2,j\}) = a$. If an edge $\{1,j\}$ is directed towards $j$ for any $j \in B'$, then $\mathcal{O}'$ has a cycle. Thus, all edges $\{1,j\}$ for $j \in B'$ are directed towards 1.

For a fixed $k \in B'$, we know $(k,1), (k,2), (k,3) \in \mathcal{O}'$. Clearly there is at least one $k' \in B \setminus \{k\}$ for which $(k',k) \in \mathcal{O}'$, and acyclicity of $\mathcal{O}'$ requires that $(k',1),(k',2) \in \mathcal{O}'$. Similarly, there must be at least one $k'' \in B \setminus \{k,k'\}$ for which $(k'',k') \in \mathcal{O}'$, with acyclicity of $\mathcal{O}'$ implying $(k'',1),(k'',2) \in \mathcal{O}'$. Continuation of this chain shows that in fact, we have $B \subseteq B'$, which is clearly impossible.

On the other hand, if more than one vertex in $A \setminus \{1,p\}$ has zero indegree when $\mathcal{O}'$ is restricted to $G_A$, then there would be a vertex in $A \setminus \{1,2,p-1,p\}$ with indegree strictly greater than $a$ under $\mathcal{O}'$, a contradiction. 
This finishes all the possibilities in  Case 3b.  

\vspace{3mm}
\fbox{\textbf{Claim 3.}} \textit{All edges between $A$-vertices and $B$-vertices have the same weight.}

\textit{Proof of Claim 3.} The same argument used to prove Claim 3 in Theorem \ref{thm:c-tree and no CV} applies, except the orientation $\mathcal{O}$ should also be such that $v$ is a source vertex in $G[B]$ under $\mathcal{O}'$. For each $i \in A$, the restriction of $\mathcal{O}_i$ to $G_B$ should be identical to the restriction of $\mathcal{O}$ to $G_B$. $\blacksquare$

\vspace{3mm}
\fbox{\textbf{Claim 4.}} \textit{If $N_G(0) \cap B \neq \emptyset$, then all $B$-vertices are adjacent to $0$ via edges of equal weight.}

\textit{Proof of Claim 4.} The same argument used to prove Claim 4 in Theorem \ref{thm:c-tree and no CV} applies. $\blacksquare$

\vspace{3mm}
Since all four statements are proven, the theorem holds. 
\end{proof}

\subsubsection{$G$ has cut vertices}
\label{subsubsection 6.2.2}

In this subsection, we complete our classification of all graphs of type $\mathscr{G}_{_\mathrm{I}}(A, B)$ by considering graphs which do have cut vertices.  

For a subset $S \subset V(G)$, let $\mathcal{CV}_G(S)$ denote the set of cut vertices of $G$ which are in $S$. Before stating the next theorem, we first prove a statement which shows that when considering cut vertices of $G$, for our purposes it is sufficient to consider vertices in the set $\mathcal{CV}_G(A \cup \{0\})$.

\begin{proposition}
\label{cut-vertices-A-and-0}
Let $G$ be a graph of type $\mathscr{G}_{_\mathrm{I}}(A, B)$ containing at least one cut vertex. Then there exists $u \in \mathcal{CV}_G(A \cup \{0\})$ such that $G \setminus \{u\}$ has a connected component $C$ with $V(C) \subseteq B$.
\end{proposition}

\begin{proof} 
 To see why $\mathcal{CV}_G(A \cup \{0\})$ must be non-empty, suppose to the contrary that all cut vertices of $G$ are in $B$. Then it follows from Corollary~\ref{cor:structure_quotient_graph} that $G / (A \cup \{0\})$ is a $c$-tree for some $c \in \mathbb{Z}^+$, and at least one connected component of $G_B$ contains a leaf vertex $j$ that is not adjacent to any vertices in $A \cup \{0\}$. Hence, $\indeg_{\mathcal{O}}(j) = c$ is a constant for every $\mathcal{O} \in \mathcal{A}(G)$. However, there must exist another $j' \in B$ which is adjacent to a vertex of $A\cup \{0\}$. In fact, $j'$ must be adjacent to more than one vertex of $A \cup\{0\}$ since $A \cup \{0\} $ contains no cut vertices of $G$. But then  there exist orientations $\mathcal{O} \neq \mathcal{O}'$ such that $\indeg_\mathcal{O}(j) \neq \indeg_{\mathcal{O}'}(j)$, a contradiction. Therefore, $G$ must contain at least one cut vertex in $A \cup \{0\}$.

Now either $G_A$ itself has cut vertices or does not, so let us examine these two cases individually.

\vspace{3mm}
\fbox{\textbf{Case 1.}} \textit{$G_A$ has cut vertices.} Then $G_A$ is an $a$-tree. In the following argument, suppose the proposition statement is false, that is, any connected component $C$ of $G_B$, where $C$ is not the singleton 0, is adjacent to at least two vertices in $A \cup \{0\}$.

Fix a connected component $C_1$ of $G_B$, where $C_1 \neq \{0\}$, and let $H_1 = G[A \cup \{0
\} \cup V(C_1)]$. We will use the fact that $H_1$ is invariant under the action of $\mathfrak{S}_p \times \mathfrak{S}_n$, where $n = |V(C_1) \cap B|$, as given by Proposition~\ref{prop:subgraph_invariance}.

If there exists a leaf vertex $v \in A$ that is not adjacent to $C_1$, then $v$ will have constant indegree under all orientations in $\mathcal{A}(H_1)$. Hence, by Lemma~\ref{lemma:invariance_maximal_GPF}, all $A$-vertices must have constant indegree in $H_1$. However, if $\{v_1,v_2\} \in A \cup \{0\}$ are two vertices adjacent to $C_1$, then $G$ must contain a cycle through $v_1$ and $v_2$ and intersecting $C_1$. In this case, whichever of $v_1,v_2$ is in $A$ will not have constant indegree under all orientations in $\mathcal{A}(H_1)$, a contradiction.

Therefore, $C_1$ is adjacent to every leaf vertex of $G_A$ (if $G_A$ has one non-zero leaf, $C_1$ must also be adjacent to $0$, else the proposition is true). This means $\mathcal{CV}_{H_1}(A \cup \{0\})$ is empty. Since this argument holds true for all connected components $C$ of $G_B$, excluding the singleton 0, we see that $\mathcal{CV}_{G}(A \cup \{0\})$ is empty, giving us the contradiction we desired.

\vspace{3mm}
\fbox{\textbf{Case 2.}} \textit{$G_A$ does not have cut vertices.} Then $G_A$ is either an $a$-cycle or the complete graph $K_{p+1}^{a,b}$. Since $\mathcal{CV}_G(A \cup \{0\})$ is non-empty, the proposition statement follows immediately. 
\end{proof}

We are now ready to present our final theorem of this subsection, which will complete our classification of all graphs of type $\mathscr{G}_{_\mathrm{I}}(A, B)$.

\begin{theorem}
\label{thm:connected and CV}
Let $G$ be a graph of type $\mathscr{G}_{_\mathrm{I}}(A, B)$. If $G$ has at least one cut vertex, then $G_A$ is an $a$-tree, an $a$-cycle, or the complete graph $K_{p+1}^{a,b}$, and one of the following holds $(a,b,c,d \in \mathbb{Z}^+)$:
\begin{enumerate}
    \item[(i)] there exists $i \in A \cup \{0\}$ such that $G[\{i\} \cup B]$ is a $c$-tree, $c$-cycle, or $K_{q+1}^{c,d}$ with root $i$, and no other vertices of $A \cup \{0\}$ are adjacent to vertices of $B$.
    \item[(ii)] there exist $i_1, i_2, \ldots, i_n \in A \cup \{0\}$ for $n > 1$ such that $G[(\bigcup_{k=1}^n \{i_k\}) \cup B]$ is a forest of $c$-trees, and no other vertices of $A \cup \{0\}$ are adjacent to vertices of $B$.
\end{enumerate}
\end{theorem}

\begin{proof} By Proposition \ref{cut-vertices-A-and-0}, we can find $u \in \mathcal{CV}_G(A \cup \{0\})$ and a connected component $C$ of $G \setminus \{u\}$ with $V(C) \subseteq B$. Consider the subgraph $H = G[\{u\} \cup V(C)]$. 
Notice that we may treat $u$ as the unique source of $H$, since there is a natural one-to-one correspondence between the set of all orientations in $\mathcal{A}(G)$ restricted to $H$ and the set of all orientations in $\mathcal{A}(H)$ with unique source $u$. 
Moreover, by Proposition~\ref{prop:subgraph_invariance} $H$ is invariant under the action of $\mathfrak{S}_{|V(C)|}$. Thus, by Theorem~\ref{thm:gaydarov and hopkins}, $H$ is a $c$-tree, a $c$-cycle, or the complete graph $K_{|V(C)|+1}^{c,d}$. As in the proof of Proposition~\ref{cut-vertices-A-and-0}, the case where $H$ is either a $c$-cycle or the complete graph $K_{|V(C)|+1}^{c,d}$ will be considered separately from the case where $H$ is a $c$-tree.

\vspace{3mm}
\fbox{\textbf{Case 1.}} \textit{$H$ is a $c$-cycle or the complete graph $K_{|V(C)|+1}^{c,d}$}. In particular, $|V(C)| \geq 2$. We wish to show that $V(C) = B$. 

Suppose to the contrary. 
Then $B-V(C) \neq \emptyset$. 
Let $i \in V(C)$, and observe that $i \in V(C) \subset \mathcal{I}(G)$. Then we can find an $\mathcal{O} \in \mathcal{A}(G)$ and a vertex $j \in B-C$ such that  both $\indeg_{\mathcal{O}}(i)$ and $\indeg_{\mathcal{O}}(j)$ are maximal, 
and thus equal (simply choose an orientation in $\mathcal{A}(G[V(G) \setminus V(C)])_{A \to (B-C)}$,  which  must have a sink $j \in B-C$. Hence $\indeg_{\mathcal{O}}(j)$ is maximal. Then extend the orientation in the obvious way into an orientation of $\mathcal{A}(G)$ in which $\indeg_{\mathcal{O}}(i)$ is also maximal). In particular, $\indeg_{\mathcal{O}}(i) = \deg_G(i)$. Now let $i' \in V(C)$ ($i' \neq i$) be adjacent to $i$. By Lemma~\ref{lemma:invariance_maximal_GPF}, there must exist $\mathcal{O}' \in \mathcal{A}(G)$ with $\indeg_{\mathcal{O}'}(i) = \deg_G(i) = \deg_G(i') = \indeg_{\mathcal{O}'}(i')$, which is impossible since $i$ and $i'$ are adjacent. Therefore, $V(C) = B$.

\vspace{3mm}
\fbox{\textbf{Case 2.}} \textit{$H$ is a $c$-tree}. Consequently, every vertex in $V(C)$ has constant indegree under all orientations in both $\mathcal{A}(H)$ and $\mathcal{A}(G)$. By Lemma~\ref{lemma:invariance_maximal_GPF}, this fact necessitates that all $B$-vertices in $G$ have constant indegree under all orientations in $\mathcal{A}(G)$, which occurs if and only if all components of $G_B$ (of which there may be more than one) are $c$-trees that are connected to a single vertex of $A \cup \{0\}$ (cf. Figure~\ref{fig:connected_cv_examples:sub2}). 
\end{proof}


\begin{figure}[H]
\begin{center}
\begin{subfigure}[t]{.45\textwidth}
\begin{center}
\begin{tikzpicture}[line width=.7pt,scale=1]
\tikzstyle{vertex}=[draw,circle,fill=black,minimum size=3pt,inner sep=0pt]
\path (0,0) node[vertex, label=left:0] (v0) {};
\path (1,1) node[vertex, label=above:1] (v1) {};
\path (2.3,1) node[vertex, label=above:2] (v2) {};
\path (2,0) node[vertex, label=left:3] (v3) {};
\path (1,-1) node[vertex, label=right:4] (v4) {};
\path (3.5,1) node[vertex, label=above:5] (v5) {};
\path (4.5,0.5) node[vertex, label=right:6] (v6) {};
\path (4.5,-0.5) node[vertex, label=right:7] (v7) {};
\path (3.5,-1) node[vertex, label=below:8] (v8) {};

\draw[blue] (v0) -- (v1) node[midway, above=1mm] {$a$}; 
\draw[blue] (v1) edge (v2) edge (v3);
\draw[blue] (v0) edge (v4);

\draw[Green] (v5) -- (v6) node[midway, above=0.5mm] {$d$};
\draw[Green] (v6) -- (v7) -- (v8) -- (v5) -- (v7) -- (v8) -- (v6);

\draw[red] (v3) -- (v5) node[midway, above=1mm] {$c$};
\draw[red] (v3) edge (v6) edge (v7) edge (v8);

\end{tikzpicture}
\caption{A graph described in Theorem \ref{thm:connected and CV}(i), where $A=\{1,2,3,4\}$ and $B=\{5,6,7,8\}$.}
\label{fig:connected_cv_examples:sub1}
\end{center}
\end{subfigure}
\hspace{5mm}
\begin{subfigure}[t]{.45\textwidth}
\begin{center}
\begin{tikzpicture}[line width=.7pt,scale=1]
\tikzstyle{vertex}=[draw,circle,fill=black,minimum size=3pt,inner sep=0pt]
\path (0,0) node[vertex, label=left:0] (v0) {};
\path (1,.5) node[vertex, label=above:1] (v1) {};
\path (2,.75) node[vertex, label=above:2] (v2) {};
\path (3,.5) node[vertex, label=above:3] (v3) {};
\path (4,0) node[vertex, label=right:4] (v4) {};
\path (.5,-1) node[vertex, label=below:5] (v5) {};
\path (1,-1) node[vertex, label=below:6] (v6) {};
\path (1.75,-1) node[vertex, label=below:7] (v7) {};
\path (2.25,-1) node[vertex, label=below:8] (v8) {};
\path (3,-1) node[vertex, label=below:9] (v9) {};

\draw[blue] (v0) -- (v1) node[midway, above=0.2mm] {$a$};
\draw[blue] (v1) -- (v2) -- (v3) -- (v4) -- (v0);

\draw[red] (v0) -- (v5) node[midway, left=0.1mm] {$c$};
\draw[red] (v5) -- (v6);
\draw[red] (v2) edge (v7) edge (v8);
\draw[red] (v3) edge (v9);
 
\end{tikzpicture}
\caption{A graph described in Theorem \ref{thm:connected and CV}(ii), where $A=\{1,2,3,4\}$ and $B=\{5,6,7,8,9\}$.}
\label{fig:connected_cv_examples:sub2}
\end{center}
\end{subfigure}
\end{center}
\caption{Examples of graphs in Theorem~\ref{thm:connected and CV}.}
\label{fig:connected_cv_examples}
\end{figure}


\subsection{Neither $G_A$ nor $G_B$ is connected}
\label{subsec:not-connect}

In this subsection, we show that if neither of the subgraphs $G_A$ nor $G_B$ is connected, then the number of permissible structures of the graph $G$ is dramatically reduced. As in the previous subsection, we begin by introducing some special notation which will be used in the proofs to follow. We will say a graph $G$ is of type $\mathscr{G}_{_\mathrm{II}}(A,B)$ if $G$ satisfies the following:
\begin{enumerate}
    \item $V(G) = \{0\} \cup A \cup B$, where 0 is the unique root vertex of $G$, $A = \{1,2,\ldots,p\}$, and $B=\{p+1,p+2,\ldots,p+q\}$ for some $p,q \in \mathbb{Z}^+$;
    \item neither $G_A$ nor $G_B$ is connected;
    \item $\cpf(G)$ is invariant under the action of the product of symmetric groups $\mathfrak{S}_{p} \times \mathfrak{S}_{p}$.
\end{enumerate}

In our classification of all graphs of type $\mathscr{G}_{_\mathrm{II}}(A, B)$, we will use a special graph $\tilde{G}$ constructed from $G$, which we will call the \textit{component graph} of $G$. Given a graph $G$, the component graph $\tilde{G}$ is the quotient graph of $G$ with respect to the equivalence relation $\sim$ defined by: $u \sim v$ if and only if, for either $S=A$ or $S=B$, $u, v \in S$ and there exists a connected component of $G_S$ containing both $u$ and $v$, with 0 in its own equivalence class. We will refer to the vertices of $\tilde{G}$ corresponding to connected components of $G_A$ (resp. $G_B$) as the \textit{$A$-nodes} (resp. \textit{$B$-nodes}) of $\tilde{G}$, and the vertex of $\tilde{G}$ containing 0 as the \textit{$0$-node}. When we say two nodes of $\tilde{G}$ (or in general, two connected components in $G_A$ or $G_B$) are adjacent, we mean there exists at least one edge between vertices in each node.
Figure~\ref{fig:component_graph_ex} shows an example of a component graph.

\begin{figure}[H]
\begin{center}
\begin{subfigure}[t]{.45\textwidth}
\begin{center}
\begin{tikzpicture}[line width=.7pt,scale=1]
\tikzstyle{vertex}=[draw,circle,fill=black,minimum size=4pt,inner sep=0pt]
\path (0,0) node[vertex, label=left:0] (v0) {};
\path (1,1.5) node[vertex, red] (v1) {};
\path (1,.5) node[vertex, red] (v2) {};
\path (2,1) node[vertex, blue] (v3) {};
\path (3,1.5) node[vertex, blue] (v4) {};
\path (3,.5) node[vertex, blue] (v5) {};
\path (4,1.5) node[vertex, red] (v6) {};
\path (4,.5) node[vertex, red] (v7) {};
\path (5,1) node[vertex, red] (v8) {};

\path (1,-1) node[vertex, blue] (v9) {};
\path (2,-1) node[vertex, blue] (v10) {};
\path (3,-1) node[vertex, blue] (v11) {};
\path (3,-.25) node[vertex, red] (v12) {};
\path (4,-.5) node[vertex, red] (v13) {};
\path (4,-1.5) node[vertex, red] (v14) {};

\draw (v0) -- (v1) -- (v3) -- (v4) -- (v6) -- (v8); 
\draw (v0) -- (v2) -- (v3) -- (v5) -- (v7);
\draw (v1) -- (v2);
\draw (v6) -- (v7);
\draw (v0) -- (v9) -- (v10) -- (v11) -- (v13) -- (v14);
\draw (v10) -- (v12);
\draw (v13) -- (v14);

\end{tikzpicture}
\caption*{}
\end{center}
\end{subfigure}
\hspace{5mm}
\begin{subfigure}[t]{.45\textwidth}
\begin{center}
\begin{tikzpicture}[line width=.7pt,scale=1]
\tikzstyle{vertex0}=[draw,circle,fill=black,minimum size=4pt,inner sep=0pt]
\tikzstyle{vertex}=[draw,circle,fill=none,minimum size=8pt,inner sep=0pt,thick]
\path (0,0) node[vertex0, label=left:0] (v0) {};
\path (1,1) node[vertex, red] (v1) {};
\path (2.5,1) node[vertex, blue] (v2) {};
\path (4,1) node[vertex, red] (v3) {};

\path (1,-1) node[vertex, blue] (v4) {};
\path (2.5,-.5) node[vertex, red] (v5) {};
\path (2.5,-1.5) node[vertex, red] (v6) {};

\draw (v0) -- (v1) -- (v2) -- (v3); 
\draw (v0) -- (v4) -- (v5);
\draw (v4) -- (v6);

\end{tikzpicture}
\caption*{}
\end{center}
\end{subfigure}
\end{center}
\caption{A graph $G$ (on left) and its component graph $\tilde{G}$ (on right).}
\label{fig:component_graph_ex}
\end{figure}

\begin{remark}
\label{remark:comp_graph_cycle_pendant}
Based on the theorems of the previous subsection, we know the possible structures of the component graphs for every graph of type $\mathscr{G}_{_\mathrm{I}}(A, B)$. In particular, if a component graph $\tilde{G}$ of a graph $G$ of type $\mathscr{G}_{_\mathrm{I}}(A, B)$ contains a cycle,
then $\tilde{G}$ cannot contain a node that is a leaf.
\end{remark}

\begin{lemma}
\label{lemma:component_graph_no_triangle}
Let $G$ be a graph of type $\mathscr{G}_{_\mathrm{II}}(A, B)$ and $\tilde{G}$ its component graph. If $\tilde{G}$ has a cycle, then $\tilde{G}$ has a cycle which goes through the $0$-node.
\end{lemma}

\begin{proof} Observe that in any cycle of $\tilde{G}$, no two $A$-nodes (resp. $B$-nodes) can be adjacent by definition of $\tilde{G}$. With this in mind, we prove the lemma by contradiction.

Suppose no cycle of $\tilde{G}$ goes through the 0-node. Then we can fix a cycle $\tilde{\mathcal{C}}$ of $\tilde{G}$ and represent it by a sequence of alternating $A$- and $B$-nodes: ($A_1, B_2, A_3, B_4, \ldots, B_{n}, A_1$). Without loss of generality, we can assume this cycle is ``closest" to the 0-node among all cycles in $\tilde{G}$ in the sense that, for any other cycle $\tilde{\mathcal{C}}'$ of $\tilde{G}$, the minimal distance between 0 and a node of $\tilde{\mathcal{C}}'$ is at least the minimal distance between 0 and a node of $\tilde{\mathcal{C}}$. We may also assume that the subgraph of $\tilde{G}$ induced by the nodes of $\tilde{\mathcal{C}}$ is a single cycle (otherwise, we may restrict to a smaller cycle in $\tilde{\mathcal{C}}$). In addition, since $\tilde{\mathcal{C}}$ does not contain the 0-node, there must be a directed path of nodes $\tilde{P} = (0, S_1, S_2, \ldots, S_m)$, $m \geq 1$, in $\tilde{G}$ from $0$ to one node $S_m$ in $\tilde{\mathcal{C}}$. 

By the structure described above, it is possible to choose vertices $v_1, v_2, \ldots, v_s \in V(G) \setminus \{0\}$ from the nodes of $\tilde{\mathcal{C}}$ (at least one vertex from each node) such that $G[\bigcup_{i=1}^s \{v_i\}]$ is a single cycle. 

\fbox{ \textbf{Case 1.}} \emph{Assume that the  path $\tilde{P}$  has length at least two}.  Then there exist $u_1, u_2, \ldots, u_t \in V(G) \setminus \{0\}$ $(t \geq 1)$ such that $G[\{0\} \cup \bigcup_{j=1}^t \{u_j\}]$ is the path $(0,u_1,u_2,\ldots,u_t)$ from $0$ to $u_t$ containing vertices from every node of $\tilde{P}$, with $u_t$ adjacent to at least one vertex in the set $\{v_1,\ldots,v_s\}$. The induced subgraph $H \coloneqq G[\{0\} \cup (\bigcup_{i=1}^s \{v_i\}) \cup (\bigcup_{j=1}^t \{u_j\})]$ will resemble Figure~\ref{fig:disconnected_comp_graph_cycle}. However, there must be at least one $u_j$ (e.g., $u_1$) with constant indegree under every orientation in $\mathcal{A}(H)$, in contrast to each $v_i$ not in the same A-component of $v_1$ having non-constant indegree. Since there are both $A$- and $B$-vertices in the cycle $(v_1, v_2, \ldots, v_s, v_1)$, this violates Proposition~\ref{prop:subgraph_invariance}.

\fbox{\textbf{Case 2.}} The path 
 $\tilde{P} = (0,S_1)$, with $S_1$ a node in $\tilde{C}$. 
 Applying Proposition~\ref{prop:subgraph_invariance} and 
 Theorem~\ref{thm:gaydarov and hopkins}, we can reduce the configuration to two subcases. 
 \begin{enumerate}[(i)]
     \item There is a path $(0, u_1, u_2, \dots, u_t)$ from $0$ to $u_t$, where $u_t$ is adjacent to $v_1$. There is no other  edges between $\{0, u_1, \dots, u_t\}$ to $\{v_1, \dots, v_s\}$. And $u_1, \dots, u_t, v_1$ are all $A$-vertices or all $B$-vertices. WOLG,  assume they are $A$-vertices. 
     This subcase is similar to Case 1.  In the induced subgraph on $H \coloneqq G[\{0\} \cup (\bigcup_{i=1}^s \{v_i\}) \cup (\bigcup_{j=1}^t \{u_j\})]$, $v_1$ has constant indegree under every orientation in $\mathcal{A}(H)$, while a $A$-vertex in other $A$-components of $\tilde{C}$  has non-constant indegree. 
     \item The root $0$ is adjacent to $v_1$, $v_2$ in the cycle, where  $v_1$ and $v_2$ are both in $A$ or both in $B$. WOLG, assume
     $v_1, v_2 \in A$. Let $H:=G[\{0\} \cup  (\bigcup_{i=1}^s \{v_i\})]$. Then in any orientation in $\mathcal{A}(H)$, $v_1$
     and $v_2$ cannot both reach maximal degree, but $v_1$ and $v_j$ can both reach the maximal degree if $v_j$ is an $A$-vertex in a different $\tilde{C}$ component than $v_1$. 
 \end{enumerate}

In any of the above cases, we reach a contradiction. 
This finishes the proof of the lemma. 
\end{proof}


\begin{figure}[H]
\begin{center}
\begin{tikzpicture}[line width=.7pt,scale=.78]
\tikzstyle{vertex0}=[draw,circle,fill=black,minimum size=4pt,inner sep=0pt]
\tikzstyle{vertex}=[draw,circle,fill=none,minimum size=8pt,inner sep=1pt,thick]
\tikzstyle{evertex}=[draw,circle,fill=none,minimum size=5pt,inner sep=2pt,thick]

\path (-4,1) node[vertex0, label=left:0] (v0) {};

\path (1.5*360/6: 2cm) node[vertex,red] (A1) {$A_1$};
\path (2.5*360/6: 2cm) node[vertex,blue] (B1) {$B_1$};
\path (3.5*360/6: 2cm) node[vertex,red] (A2) {$A_2$};
\path (4.5*360/6: 2cm) node[vertex,blue] (B2) {$B_2$};
\path (5.5*360/6: 2cm) node[vertex,red] (A3) {$A_3$};
\path (.5*360/6: 2cm) node[vertex,blue] (B3) {$B_3$};

\path (-3,2) node[vertex,red] (S1) {$S_1$};
\path (-1.5,2.5) node[vertex,blue] (S2) {$S_2$};

\draw (A1) edge[bend right=15] (B1);
\draw (B1) edge[bend right=15] (A2);
\draw (A2) edge[bend right=15] (B2);
\draw (B2) edge[bend right=15] (A3);
\draw (A3) edge[bend right=15] (B3);
\draw (B3) edge[bend right=15] (A1);
\draw (v0) edge[bend left=15] (S1);
\draw (S1) edge[bend left=15] (S2);
\draw (S2) edge[bend left=15] (A1);

\draw[->] (3.5,0) -- (5.5,0);

\path (6,1) node[vertex0,label=left:0] (v00) {};
\path (7,2) node[vertex0,red,label=above:$u_1$] (S11) {};
\path (8.3,2.5) node[vertex0,blue,label=above:$u_2$] (S22) {};
\path (8.9,2.5) node[vertex0,blue,label=above:$u_3$] (S23) {};

\path (1.3*360/6: 2cm)+(10,0) node[vertex0,red,label=below:$v_1$] (v1) {};
\path (1.7*360/6: 2cm)+(10,0) node[vertex0,red,label=below:$v_2$] (v2) {};
\path (2.5*360/6: 2cm)+(10,0) node[vertex0,blue,label=left:$v_3$] (v3) {};
\path (3.3*360/6: 2cm)+(10,0) node[vertex0,red,label=left:$v_4$] (v4) {};
\path (3.7*360/6: 2cm)+(10,0) node[vertex0,red,label=left:$v_5$] (v5) {};
\path (4.1*360/6: 2cm)+(10,0) node[vertex0,blue,label=below left:$v_6$] (v6) {};
\path (4.5*360/6: 2cm)+(10,0) node[vertex0,blue,label=below:$v_7$] (v7) {};
\path (4.9*360/6: 2cm)+(10,0) node[vertex0,blue,label=below right:$v_8$] (v8) {};
\path (5.5*360/6: 2cm)+(10,0) node[vertex0,red,label=right:$v_9$] (v9) {};
\path (.3*360/6: 2cm)+(10,0) node[vertex0,blue,label=right:$v_{10}$] (v10) {};
\path (.7*360/6: 2cm)+(10,0) node[vertex0,blue,label=above right:$v_{11}$] (v11) {};

\draw (v1) edge[bend right=15] (v2);
\draw (v2) edge[bend right=15] (v3);
\draw (v3) edge[bend right=15] (v4);
\draw (v4) edge[bend right=15] (v5);
\draw (v5) edge[bend right=15] (v6);
\draw (v6) edge[bend right=15] (v7);
\draw (v7) edge[bend right=15] (v8);
\draw (v8) edge[bend right=15] (v9);
\draw (v9) edge[bend right=15] (v10);
\draw (v10) edge[bend right=15] (v11);
\draw (v11) edge[bend right=15] (v1);
\draw (v00) edge[bend left=15] (S11);
\draw (S11) edge[bend left=15] (S22);
\draw (S22) edge[bend left=15] (S23);
\draw (S23) edge[bend left=15] (v1);
\draw[dashed] (S23) edge[bend right=15] (v2);

\end{tikzpicture}
\end{center}
\caption{An example of the cycle $\tilde{C}$ and path $\tilde{P}$ (on left), and a subgraph of $G$ induced by vertices in nodes of $\tilde{C}$ and $\tilde{P}$ (on right).}
\label{fig:disconnected_comp_graph_cycle}
\end{figure}


\begin{theorem}
\label{thm:disconnected_comp_graph_cycle}
Let $G$ be a graph of type $\mathscr{G}_{_\mathrm{II}}(A, B)$ and $\tilde{G}$ its component graph. If $\tilde{G}$ has a cycle, then $G$ is an $a$-cycle for some $a \in \mathbb{Z}^+$.
\end{theorem}

\begin{proof} 
By Lemma~\ref{lemma:component_graph_no_triangle}, $\tilde{G}$ has a cycle $\tilde{\mathcal{C}}$ which goes through 0, which can be represented by the sequence of alternating $A$- and $B$-nodes $(0, S_1, S_2, \ldots, S_n, 0)$. As in the proof of Lemma~\ref{lemma:component_graph_no_triangle}, we can assume the subgraph of $\tilde{G}$ induced by the nodes of $\tilde{\mathcal{C}}$ is a single cycle, and we can again find vertices $v_1, v_2, \ldots, v_s \in V(G)$ such that $H \coloneqq G[\{0\} \cup (\bigcup_{i=1}^s \{v_i\})]$ is a single cycle containing vertices from every node of $\tilde{\mathcal{C}}$. We will first show that $\tilde{C}$ must have minimum length five.

\vspace{3mm}
\fbox{\textbf{Claim 1.}} \textit{$\tilde{C}$ must have length greater than three.} 

\textit{Proof of Claim 1.} Suppose $\tilde{C}$ has length three, with $\tilde{C} = (0, A_1, B_1, 0)$. 
Since $G_B$ is disconnected, there must exist another $B$-node $B_2$ in $\tilde{G}$ such that $B_2$ is not adjacent to the 0-node, and $B_2$ satisfies one of the following: either $B_2$ is adjacent to $A_1$ (cf. Figure~\ref{fig:proof_comp_graph_cycle_1a}), or else there is another $A$-node $A_2$ which is adjacent to $B_2$ and the $0$-node (cf. Figure~\ref{fig:proof_comp_graph_cycle_1b}). Define $H$ to be the subgraph $G$ induced by vertices in the $0$-node, $A_1$, $B_1$, $B_2$, and (if required) $A_2$. Then $H$ is of type $\mathscr{G}_{_\mathrm{I}}(A^*, B^*)$, where $A* = V(H) \cap A$ and $B^* = V(H) \cap B$. Therefore, the structure of $H$ is given by one of Theorems~\ref{thm:c-cycle and no CV}-\ref{thm:complete and no CV} and 
\ref{thm:connected and CV}.

Clearly, Theorem~\ref{thm:c-cycle and no CV} does not apply, nor can Theorems~\ref{thm:c-tree and no CV}-\ref{thm:complete and no CV} apply since some $B$-vertices in $H$ are adjacent to $0$ while others are not. Likewise, Theorem~\ref{thm:connected and CV} cannot apply since $H$ has an induced cycle containing 0, at least one $A$-vertex, and at least one $B$-vertex. We conclude that $\tilde{C}$ must have length greater than three.  
$\blacksquare$

\begin{figure}[H]
\begin{center}
\begin{subfigure}[t]{.45\textwidth}
\begin{center}
\begin{tikzpicture}[line width=.7pt,scale=.75]
\tikzstyle{vertex0}=[draw,circle,fill=black,minimum size=4pt,inner sep=0pt]
\tikzstyle{vertex}=[draw,circle,fill=none,minimum size=8pt,inner sep=1pt,thick]

\path (0,0) node[vertex0, label=above:0] (v0) {};

\path (-1,-1) node[vertex,red] (A1) {$A_1$};
\path (1,-1) node[vertex,blue] (B1) {$B_1$};
\path (-1,-2.5) node[vertex,blue] (B2) {$B_2$};

\draw (v0) edge (B1);
\draw (v0) edge (A1);
\draw (A1) edge (B1);
\draw (A1) edge (B2);

\end{tikzpicture}
\caption{$B_2$ is adjacent to $A_1$.}
\label{fig:proof_comp_graph_cycle_1a}
\end{center}
\end{subfigure}
\hspace{5mm}
\begin{subfigure}[t]{.45\textwidth}
\begin{center}
\begin{tikzpicture}[line width=.7pt,scale=.75]
\tikzstyle{vertex0}=[draw,circle,fill=black,minimum size=4pt,inner sep=0pt]
\tikzstyle{vertex}=[draw,circle,fill=none,minimum size=8pt,inner sep=1pt,thick]

\path (0,0) node[vertex0, label=above:0] (v0) {};

\path (0,-1) node[vertex,red] (A1) {$A_1$};
\path (1.5,-1) node[vertex,blue] (B1) {$B_1$};
\path (-1.5,-1) node[vertex,red] (A2) {$A_2$};
\path (-1.5,-2.5) node[vertex,blue] (B2) {$B_2$};

\draw (v0) edge (B1);
\draw (v0) edge (A1);
\draw (v0) edge (A2);
\draw (A1) edge (B1);
\draw (A2) edge (B2);
\draw[dashed] (A2) edge[bend right=45] (B1);

\end{tikzpicture}
\caption{$B_2$ is adjacent to $A_2$, with $A_2$ also adjacent to the $0$-node and possibly $B_1$.}
\label{fig:proof_comp_graph_cycle_1b}
\end{center}
\end{subfigure}
\end{center}
\caption{The two ways $B_2$ can be connected to the cycle $\tilde{C}$.}
\label{fig:proof_comp_graph_cycle_1}
\end{figure}

\fbox{\textbf{Claim 2.}} \textit{$\tilde{C}$ must have length greater than four.}

\textit{Proof of Claim 2.} If $\tilde{C}$ has length four, then there exist two $A$-nodes $A_1, A_2$ and one $B$-node $B_1$ in $\tilde{G}$ such that $\tilde{C}$ is the cycle $(0,A_1,B_1,A_2,0)$ (the case with two $B$-nodes will follow by symmetry). Since $G_A$ is disconnected, there must exist another $A$-node $A_3$ which is not adjacent to the 0-node in $\tilde{G}$ and which satisfies one of the following cases:

\indent \indent \fbox{\textit{Case 2a.}} \textit{$A_3$ is adjacent to $B_1$} (cf. Figure~\ref{fig:proof_comp_graph_cycle_2a}). We know that the maximal indegree obtainable by all $A$-vertices in $G[V(H) \cup V(A_3)]$ under orientations in $\mathcal{A}(G[V(H) \cup V(A_3)])$ must be the same number, say $\alpha$. Observe that if $u \in V(A_1)$ and $v \in V(A_3)$, then it is clear from the structure of $G[V(H) \cup V(A_3)]$ that there exists an $\mathcal{O} \in \mathcal{A}(G[V(H) \cup V(A_3)])$ where $\indeg_{\mathcal{O}}(u) = \indeg_{\mathcal{O}}(v) = \alpha$. However, if $u' \in V(A_2)$, there is no $\mathcal{O}' \in \mathcal{A}(G[V(H) \cup V(A_3)])$ under which $\indeg_{\mathcal{O}'}(u) = \indeg_{\mathcal{O}}(u) = \alpha$ and $\indeg_{\mathcal{O}'}(u') = \indeg_{\mathcal{O}}(v) = \alpha$. Thus, this case is impossible.

\indent \indent \fbox{\textit{Case 2b.}} \textit{There is another $B$-node $B_2$ that is adjacent to $A_3$ and at least one of $0$, $A_1$, and $A_2$} (cf. Figure~\ref{fig:proof_comp_graph_cycle_2b}). Note that $B_2$ cannot be adjacent to more than one of $0$, $A_1$, and $A_2$, else this case restricts to Case 2a or $\tilde{G}$ contains a length-three cycle. Define $\tilde{H}$ to be the subgraph of $\tilde{G}$ induced by the $0$-node, $A_1$, $A_2$, $B_1$, and $B_2$. Then the subgraph of $G$ corresponding to $\tilde{H}$ is of type $\mathscr{G}_{_\mathrm{I}}(A^*, B^*)$, where $A* = V(\tilde{H}) \cap A$ and $B^* = V(\tilde{H}) \cap B$; but $\tilde{H}$ contains a cycle and the leaf node $B_2$, which is not possible per Remark~\ref{remark:comp_graph_cycle_pendant}.

\indent \indent \fbox{\textit{Case 2c.}} \textit{There is another $A$-node $A_4$ that is adjacent to 0 and another $B$-node $B_2$ that is adjacent to $A_4$ and $A_3$ but not 0} (cf. Figure~\ref{fig:proof_comp_graph_cycle_2c}). Observe that if $B_2$ is adjacent to either $A_1$ or $A_2$, this case reduces to Case 2b. So we may assume $B_2$ is only adjacent to $A_3$ and $A_4$. Mirroring the argument in Case 2b, let $\tilde{H}$ be the subgraph of $\tilde{G}$ induced by the $0$-node, $A_1$, $A_2$, $A_4$, $B_1$, and $B_2$. Then the subgraph of $G$ corresponding to $\tilde{H}$ is of type $\mathscr{G}_{_\mathrm{I}}(A^*, B^*)$, where $A* = V(\tilde{H}) \cap A$ and $B^* = V(\tilde{H}) \cap B$; yet $\tilde{H}$ contains a cycle and the leaf node $B_2$, which renders this case impossible (see Remark~\ref{remark:comp_graph_cycle_pendant}). 

Claim  2 follows by Cases 2a-2c. $\blacksquare$

\begin{figure}[H]
\begin{center}
\begin{subfigure}[t]{.3\textwidth}
\begin{center}
\begin{tikzpicture}[line width=.7pt,scale=.75]
\tikzstyle{vertex0}=[draw,circle,fill=black,minimum size=4pt,inner sep=0pt]
\tikzstyle{vertex}=[draw,circle,fill=none,minimum size=8pt,inner sep=1pt,thick]

\path (0,0) node[vertex0, label=above:0] (v0) {};

\path (-1.25,-1) node[vertex,red] (A1) {$A_1$};
\path (1.25,-1) node[vertex,red] (A2) {$A_2$};
\path (0,-2) node[vertex,blue] (B1) {$B_1$};
\path (1.25,-3) node[vertex,red] (A3) {$A_3$};

\draw (v0) edge (A1);
\draw (v0) edge (A2);
\draw (A1) edge (B1);
\draw (A2) edge (B1);
\draw (B1) edge (A3);

\end{tikzpicture}
\caption{$A_3$ is adjacent to $B_1$.}
\label{fig:proof_comp_graph_cycle_2a}
\end{center}
\end{subfigure}
\hspace{5mm}
\begin{subfigure}[t]{.3\textwidth}
\begin{center}
\begin{tikzpicture}[line width=.7pt,scale=.75]
\tikzstyle{vertex0}=[draw,circle,fill=black,minimum size=4pt,inner sep=0pt]
\tikzstyle{vertex}=[draw,circle,fill=none,minimum size=8pt,inner sep=1pt,thick]

\path (0,0) node[vertex0, label=above:0] (v0) {};

\path (0,-1.25) node[vertex,red] (A1) {$A_1$};
\path (2,-1.25) node[vertex,red] (A2) {$A_2$};
\path (1,-2.5) node[vertex,blue] (B1) {$B_1$};
\path (-2,-1.25) node[vertex,blue] (B2) {$B_2$};
\path (-2,-2.5) node[vertex,red] (A3) {$A_3$};

\draw (v0) edge (A1);
\draw (v0) edge (A2);
\draw[dashed] (v0) edge (B2);
\draw (A1) edge (B1);
\draw (A2) edge (B1);
\draw (B2) edge (A3);
\draw[dashed] (B2) edge (A1);
\draw[dashed] (B2) edge[bend right=30] (A2);

\end{tikzpicture}
\caption{$A_3$ is adjacent to $B_2$, which is adjacent to at least one of $0$, $A_1$, and $A_2$.}
\label{fig:proof_comp_graph_cycle_2b}
\end{center}
\end{subfigure}
\hspace{5mm}
\begin{subfigure}[t]{.3\textwidth}
\begin{center}
\begin{tikzpicture}[line width=.7pt,scale=.75]
\tikzstyle{vertex0}=[draw,circle,fill=black,minimum size=4pt,inner sep=0pt]
\tikzstyle{vertex}=[draw,circle,fill=none,minimum size=8pt,inner sep=1pt,thick]

\path (0,0) node[vertex0, label=above:0] (v0) {};

\path (0,-1.25) node[vertex,red] (A1) {$A_1$};
\path (2,-1.25) node[vertex,red] (A2) {$A_2$};
\path (1,-2.5) node[vertex,blue] (B1) {$B_1$};
\path (-2,-1.25) node[vertex,red] (A4) {$A_4$};
\path (-2,-2.5) node[vertex,blue] (B2) {$B_2$};
\path (-2,-3.75) node[vertex,red] (A3) {$A_3$};

\draw (v0) edge (A1);
\draw (v0) edge (A2);
\draw (v0) edge (A4);
\draw (A1) edge (B1);
\draw (A2) edge (B1);
\draw (A4) edge (B2);
\draw (B2) edge (A3);
\draw[dashed] (B2) edge (A1);
\draw[dashed] (B2) edge[bend right=5] (A2);
\draw[dashed] (A4) edge (B1);

\end{tikzpicture}
\caption{$A_3$ is adjacent to $B_2$, which is not adjacent to $0$.}
\label{fig:proof_comp_graph_cycle_2c}
\end{center}
\end{subfigure}
\end{center}
\caption{The three ways $A_3$ can be connected to the cycle $\tilde{C}$.}
\label{fig:proof_comp_graph_cycle_2}
\end{figure}

Now we will show that $G=H$ via contradiction. If $G \neq H$, then there is $x \in V(G) \setminus V(H)$ that is adjacent to at least one vertex of $H$. As noted in Claim 2, if $x$ has only one neighbor in $H$, then $x$ will have constant indegree in the induced subgraph $G[V(H) \cup \{x\}]$ under all orientations in $\mathcal{A}(G[V(H) \cup \{x\}])$, violating Proposition~\ref{prop:subgraph_invariance}. So $x$ must have at least two neighbors in $H$. 

Without loss of generality, say $x \in A$, and let $S_A$ (resp. $S_B$) be the set of $A$-vertices (resp. $B$-vertices) in $H$ adjacent to $x$. Our goal is to show that both $S_A$ and $S_B$ must be empty.

\vspace{3mm}
\fbox{\textbf{Claim 3.}} \textit{$S_A$ is empty.}

\textit{Proof of Claim 3.}  
We observe  $S_A \neq V(H) \cap A$, otherwise $\tilde{C}$ only has one $A$-node and cannot have length $\geq 5$. Hence 
we can assume $|S_A| \geq 1$ but $S_A \neq V(H) \cap A$. 
Then there exist vertices $y \in S_A$ and $y' \in (V(H) \cap A) \setminus S_A$ with $deg_H(y) = deg_H(y')$ since $y,y' \in \mathcal{I}(H) \cap A$.
However, the degree of $y$ in $G[V(H) \cup \{x\}]$ must be strictly greater than the degree of $y'$ in $G[V(H) \cup \{x\}]$. 
It would violate Lemma~\ref{lemma:non_cut_vertex_equal_degrees} 
Since $x$ has at least two neighbors in $H$, and hence, $G[V(H) \cup \{x\}]$ has no cut vertices. 

\vspace{3mm}
\fbox{\textbf{Claim 4.}} \textit{$S_B$ is empty.}

\textit{Proof of Claim 4.}

The same argument used in Claim 3 to show $S_A \neq  V(H) \cap A$ and $|S_A| \geq 1$ also shows the case that $S_B \neq  V(H) \cap B$ and $|S_B| \geq 1$ (simply replace $A$ by $B$ in the argument). 

Next we assume $S_B=V(H) \cap B$. 
Note that if $x$ is adjacent to $0$ and $A_x$ is the $A$-node of $\tilde{G}$ containing $x$, $\tilde{C}$ contains a smaller cycle of length three or four through the $0$-node. Thus, by Claim 1 and Claim 2, $x$ cannot be adjacent to $0$ when $|S_B| \geq 1$.

Now we consider the possible pairs of nodes adjacent to the $0$-node in $\tilde{C}$.

\indent \indent \fbox{\textit{Case 4a.}} \textit{The $0$-node is adjacent to two $B$-nodes $B_1, B_2$}. Then $\tilde{C}$ contains the cycle $(0,B_1,A_x,B_2,0)$ of length four. By Claim 2, this case cannot happen.

\indent \indent \fbox{\textit{Case 4b.}} \textit{The $0$-node is adjacent to two $A$-nodes $A_1, A_2$} (cf. Figure~\ref{fig:proof_comp_graph_cycle_4b}). Let $x' \in V(A_2)$ be adjacent to 0 in $H$. Then in the induced subgraph $G[(V(H)\cup\{x\}) \setminus V(A_1)]$, $x'$ will have constant indegree under all orientations in $\mathcal{A}(G[(V(H)\cup\{x\})  \setminus V(A_1)])$. However, if $x'' \in V(H) \cap A$ is not in $A_1$ or $A_2$, $x''$ will not have constant indegree under the same orientations since it is part of a cycle in $G[(V(H)\cup\{x\}) \setminus V(A_1)]$. By Proposition~\ref{prop:subgraph_invariance}, this case is impossible.

\indent \indent \fbox{\textit{Case 4c.}} \textit{The $0$-node is adjacent to an $A$-node $A_1$ and a $B$-node $B_1$} (cf. Figure~\ref{fig:proof_comp_graph_cycle_4c}). Let $y \in V(B_1)$ be adjacent to 0 in $H$. Specifically, in $G[(V(H)\cup\{x\}) \setminus V(A_1)]$, the vertex $y$ will have constant indegree under all orientations in $\mathcal{A}(G[(V(H)\cup\{x\}) \setminus V(A_1)])$. Yet there is some $y' \in V(H) \cap B$ that is not in $B_1$ which will not have constant indegree under the same orientations due to its membership in an induced cycle of $G[(V(H)\cup\{x\}) \setminus V(A_1)]$.

In all of the above cases, our contradiction stems from the assumption that $|S_B| \geq 1$. $\blacksquare$


\begin{figure}[H]
\begin{center}
\begin{subfigure}[t]{.45\textwidth}
\begin{center}
\begin{tikzpicture}[line width=.7pt,scale=.77]
\tikzstyle{vertex0}=[draw,circle,fill=black,minimum size=4pt,inner sep=0pt]
\tikzstyle{vertex}=[draw,circle,fill=none,minimum size=35pt,inner sep=1pt,thick]

\path (1.5*360/6: 2cm) node[vertex0, label=above:0] (v0) {};
\path (0,0) node[vertex0, label=above:$x$] (x) {};

\path (2.4*360/6: 2cm) node[vertex0,red] (A11) {};
\path (2.6*360/6: 2cm) node[vertex0,red] (A12) {};
\path (3.3*360/6: 2cm) node[vertex0,blue] (B11) {};
\path (3.5*360/6: 2cm) node[vertex0,blue] (B12) {};
\path (3.7*360/6: 2cm) node[vertex0,blue] (B13) {};
\path (4.5*360/6: 2cm) node[vertex0,red,label=above:$x''$] (x'') {};
\path (5.4*360/6: 2cm) node[vertex0,blue] (B21) {};
\path (5.6*360/6: 2cm) node[vertex0,blue] (B22) {};
\path (6.5*360/6: 2cm) node[vertex0,red,label=above:$x'$] (x') {};

\path (2.5*360/6: 2cm) node[vertex,red,label=left:\textcolor{red}{$A_1$}] (A1) {};
\path (3.5*360/6: 2cm) node[vertex,blue,label=left:\textcolor{blue}{$B_1$}] (B1) {};
\path (4.5*360/6: 2cm) node[vertex,red,label=below:\textcolor{red}{$A_3$}] (A3) {};
\path (5.5*360/6: 2cm) node[vertex,blue,label=right:\textcolor{blue}{$B_2$}] (B2) {};
\path (6.5*360/6: 2cm) node[vertex,red,label=right:\textcolor{red}{$A_2$}] (A2) {};

\draw (v0) edge[bend right=15] (A11); \draw (A11) edge[bend right=15] (A12); \draw (A12) edge[bend right=15] (B11); \draw (B11) edge[bend right=15] (B12); \draw (B12) edge[bend right=15] (B13); \draw (B13) edge[bend right=15] (x''); \draw (x'') edge[bend right=15] (B21); \draw (B21) edge[bend right=15] (B22); \draw (B22) edge[bend right=15] (x'); \draw (x') edge[bend right=15] (v0);

\draw (x) edge (B11) edge (B12) edge (B13) edge (B21) edge (B22);

\end{tikzpicture}
\caption{The $0$-node is adjacent to two $A$-nodes, as in Case 4b.}
\label{fig:proof_comp_graph_cycle_4b}
\end{center}
\end{subfigure}
\hspace{5mm}
\begin{subfigure}[t]{.45\textwidth}
\begin{center}
\begin{tikzpicture}[line width=.7pt,scale=.77]
\tikzstyle{vertex0}=[draw,circle,fill=black,minimum size=4pt,inner sep=0pt]
\tikzstyle{vertex}=[draw,circle,fill=none,minimum size=35pt,inner sep=1pt,thick]

\path (1.25*360/5: 2cm) node[vertex0, label=above:0] (v0) {};
\path (0,0) node[vertex0, label=above:$x$] (x) {};

\path (2.15*360/5: 2cm) node[vertex0,red] (A11) {};
\path (2.35*360/5: 2cm) node[vertex0,red] (A12) {};
\path (3.25*360/5: 2cm) node[vertex0,blue,label=left:$y'$] (y') {};
\path (4.25*360/5: 2cm) node[vertex0,red] (x'') {};
\path (5.1*360/5: 2cm) node[vertex0,blue] (B21) {};
\path (5.35*360/5: 2cm) node[vertex0,blue,label=right:$y$] (y) {};

\path (2.25*360/5: 2cm) node[vertex,red,label=left:\textcolor{red}{$A_1$}] (A1) {};
\path (3.25*360/5: 2cm) node[vertex,blue,label=below:\textcolor{blue}{$B_2$}] (B2) {};
\path (4.25*360/5: 2cm) node[vertex,red,label=below:\textcolor{red}{$A_3$}] (A3) {};
\path (5.25*360/5: 2cm) node[vertex,blue,label=right:\textcolor{blue}{$B_1$}] (B1) {};

\draw (v0) edge[bend right=15] (A11); \draw (A11) edge[bend right=15] (A12); \draw (A12) edge[bend right=15] (y'); \draw (y') edge[bend right=15] (x''); \draw (x'') edge[bend right=15] (B21); \draw (B21) edge[bend right=15] (y); \draw (y) edge[bend right=15] (v0);

\draw (x) edge (y') edge (B21) edge (y);

\end{tikzpicture}
\caption{The $0$-node is adjacent to an $A$-node and a $B$-node, as in Case 4c.}
\label{fig:proof_comp_graph_cycle_4c}
\end{center}
\end{subfigure}
\end{center}
\caption{Examples of the induced subgraph $G[V(H)\cup\{x\}]$ in Case 4b and Case 4c.}
\label{fig:proof_comp_graph_cycle_4}
\end{figure}


Therefore, we conclude $G = H$. It remains to show that $G$ has all equally weighted edges. Clearly, if not all edges of $G$ are of equal weight, there can only be two possible weights $a,b \in \mathbb{Z}^+$ such that every vertex of $G$ is adjacent to one edge of weight $a$ and one edge of weight $b$. In particular, for any vertex $v \in V(G)$ not adjacent to 0, there are three possible indegrees for $v$ under orientations in $\mathcal{A}(G)$ when $a \neq b$: $a$, $b$, and $a+b$.  However, if $v$ is adjacent to 0, there are only two possible indegrees for $v$ under orientations in $\mathcal{A}(G)$: $a+b$ and the weight of the edge it shares with 0. It forces $a=b$ and so $G$ is an $a$-cycle. 
\end{proof}

\begin{theorem}
\label{thm:disconnected_comp_graph_tree}
Let $G$ be a graph of type $\mathscr{G}_{_\mathrm{II}}(A, B)$ and $\tilde{G}$ its component graph. If $\tilde{G}$ is a tree, then one of the following holds ($a,b \in \mathbb{Z}^+$):
\begin{enumerate}
    \item[(i)] $G$ is an $(a,b)$-tree; or
    \item[(ii)] $\tilde{G}$ has only one $B$-node, and there exists $i \in A$ such that $G[\{i\} \cup B]$ is a $b$-cycle or the complete graph $K_{q+1}^{b,c}$ with root $i$ (or $A$ and $B$ are reversed). All other nodes of $\tilde{G}$ correspond to $a$-trees. (cf. Figure~\ref{fig:case-6}.)
\end{enumerate}
\end{theorem}

\begin{proof} Without loss of generality, assume at least one $A$-vertex is adjacent to $0$ in $G$. We first
show a lemma. 

\begin{lemma}  \label{lem:A-vertex-tree}
Under the condition described above, 
any non-zero node of $\tilde{G}$ that is adjacent to the $0$-node cannot be adjacent to $0$ by more than one edge in $G$. 
\end{lemma} 
\noindent \textit{Proof of Lemma~\ref{lem:A-vertex-tree}.} \ 
Suppose the $A$-node $A_1$ is adjacent to the $0$-node in $\tilde{G}$, and assume there are two vertices $u,u'$ in $A_1$ which are adjacent to $0$ in $G$.
Then, the subgraph $H_{A_1}$ of $G$ induced by $0$ and $A_1$ must be either an $a$-cycle or the complete graph $K_{|V(A_1)|+1}^{a,b}$, and there is a vertex $u''$ in $A_1$ that adjacent to $u$. 
In any $\mathcal{O} \in \mathcal{A}(G)$, 
$\indeg_{\mathcal{O}}(u)$ and $\indeg_{\mathcal{O}}(u'')$ cannot both reach maximal under $\mathcal{O}$. 
As $G_A$ is disconnected, there must be another $v \in A$ which is not in $A_1$. By the tree structure of $\tilde{G}$, there should be an $\mathcal{O}' \in \mathcal{A}(G)$ in which both $\indeg_{\mathcal{O}'}(u)$ and $\indeg_{\mathcal{O}'}(v)$ are maximal, and thus equal. This is a violation of the symmetry of 
Lemma~\ref{lemma:invariance_maximal_GPF}. \ $\blacksquare$

If $u$ is the vertex of $A_1$ adjacent to $0$ in $G$, then $u$ will have constant indegree among all orientations in $\mathcal{A}(G)$. Hence, all $A$-vertices of $G$ must have constant indegree $a$, which necessitates that all $A$-nodes in $\tilde{G}$ correspond to $a$-trees in $G$.

We will now further analyze the structure of $G$ according to whether $\tilde{G}$ has exactly one $B$-node or more than one $B$-node.

\vspace{3mm}
\fbox{\textbf{Case 1.}} \textit{$\tilde{G}$ has exactly one $B$-node $B_1$.} Since $G_B$ is assumed to be disconnected, there must exist one $A$-node $A_1$ in $\tilde{G}$ that is adjacent to both the $0$-node and $B_1$. Moreover, there can only be one vertex $u$ of $A_1$ that is adjacent to vertices in $B_1$ due to all $A$-vertices in $G$ having constant indegree under all orientations in $\mathcal{A}(G)$. Importantly, if we restrict to the subgraph $H = G[\{u\} \cup B]$ of $G$ (note that $B_1$ contains all $B$-vertices), we see that $u$ may be treated as the unique source of $H$ due to $u$ being a cut vertex. Thus, by Proposition~\ref{prop:subgraph_invariance}, $H$ is invariant under the action of $\mathfrak{S}_{|B|}=\mathfrak{S}_q$, so that $H$ is a $b$-tree or $b$-cycle or the complete graph $K_{q+1}^{b,c}$. 

Additionally, since $G_A$ is also disconnected, there must exist at least one more $A$-node, distinct from $A_1$, that is adjacent to $B_1$ in $\tilde{G}$. Since all $A$-vertices in $G$ has constant indegree $a$, we see that any such $A$-node must be adjacent to $B_1$ via exactly one edge. See Figure~\ref{fig:case-6} for an example of $G$.

\vspace{3mm}
\fbox{\textbf{Case 2.}} \textit{$\tilde{G}$ has more than one $B$-node.}
Since $\tilde{G}$ is a tree, by the proof of Lemma~\ref{lem:A-vertex-tree}, each B-node $B_i$ must be connected a unique vertex $u_i \in A\cup \{0\}$,  where $u_i$ lies between $0$ and $B_i$ and can be viewed as the root for $B_i$.  Each induced graph on $B_i \cup \{u_i\}$  must be a $b$-tree or $b$-cycle or the complete graph $K^{b,c}_{t+1}$ for some $t \geq 2$. 

If there exists an $i$ such that $B_i \cup \{u_i\}$ is a $b$-cycle or 
$K^{b,c}_{t+1}$, then similar to the proof of Lemma~\ref{lem:A-vertex-tree}, there are two adjacent vertices $u, u''$ in $B_i$ that cannot both reach the maximal indegree in any orientation of $\mathcal{A}(G)$, yet there is a $B$-vertex $v \not \in B_i$  that can reach maximal indegree with $u$ in some orientation $\mathcal{O'} \in \mathcal{A}(G)$. This gives a contradiction. 

Thus, every $B_i \cup \{u_i\}$ is a $b$-tree. 
We conclude that $G$ must be an $(a,b)$-tree. 
\end{proof}

\section{Overlap between $G$-parking functions and  $\bsy{U}$-parking functions}
\label{sec:finish}

Having completed our classification of all graphs $G$ whose  set of $G$-parking functions are invariant under $\mathfrak{S}_p \times \mathfrak{S}_q$, it remains to show that each such graph $G$ corresponds to  a 2-dimensional weight set $\bsy{U}$ for which $\cpf_{p,q}^{(2)}(\bsy{U}) = \cpf(G)$. Our next  theorem establishes this correspondence. 
For each case, we present a representative example of  $\bsy{U}$. We then provide a complete characterization of all possible  $\bsy{U}$ by analyzing the set $\mathcal{MI}_{p,q}^{(2)}(\bsy{U})$ of increasing maximal $\bsy{U}$-parking functions.

\begin{theorem}
\label{maintheorem2}
Let $G$ be a graph with $V(G) = \{0\} \cup A \cup B$, where $0$ is the unique root vertex of $G$, $A=\{1, \dots, p\}$, and  $B=\{p+1, \dots, p+q\}$ for $p,q \in \mathbb{Z}^+$. Assume $N_G(0) \cap A \neq \emptyset$.
If $\cpf(G)$ is invariant under the action of $\mathfrak{S}_p \times \mathfrak{S}_q$, then one of the following holds:
\begin{enumerate}
    \item[(i)] $G$ is described by Theorem~\ref{maintheorem}(i) and $\cpf(G) = \cpf_{p,q}^{(2)}(\bsy{U})$, where $\bsy{U}$
     is given by:
    \begin{enumerate}
        \item[(a)]  if $p \leq 2$, then $u_{p-1,q} = a+b$, $v_{p,q-1} = 2b$, $u_{i,j} = a$ for $(i,j) \neq (p-1,q)$, and $v_{i,j} = b$ for $(i,j) \neq (p,q-1)$ ; or
        \item[(b)] if $p > 2$, then $u_{p-1,q} = v_{p,q-1} = 2a$ and $u_{i,j} = v_{i,j} = a$ for all other $(i,j)$.
    \end{enumerate}

    \item[(ii)] $G$ is described by Theorem~\ref{maintheorem}(ii) and  $\cpf(G) = \cpf_{p,q}^{(2)}(\bsy{U})$, where $p=2$, and $\bsy{U} $
    is given by 
    \[\left\{ 
    \begin{array}{ll} 
    u_{0,j} = a   & \text{ for }  j \leq q, \\
    u_{1,j} = a+b &  \text{ for } j \leq q-1, \\ 
    u_{1,q}=a+b+c, &   \\ 
    v_{i,j}=c & \text{ for }  (i,j) \neq (2, q-1), \\  
    v_{2,q-1} = 2c. & 
    \end{array} 
    \right. 
    \]

    \item[(iii)] $G$ is described by Theorem ~\ref{maintheorem}(iii), and $\cpf(G) = \cpf_{p,q}^{(2)}(\bsy{U})$ where $\bsy{U} $ 
   is given by 
   \begin{eqnarray}  \label{Eq:affine-matrix}
   \begin{pmatrix}
       u_{i,j} \\ v_{i,j} 
   \end{pmatrix}
   = \begin{pmatrix}
     b & c \\ c & d 
     \end{pmatrix}
     \begin{pmatrix}
       i \\ j
   \end{pmatrix}
  +  \begin{pmatrix}
       a \\ e
   \end{pmatrix}. 
   \end{eqnarray}

    \item[(iv)] $G$ is described by Theorem~\ref{maintheorem}(iv) and $\cpf(G) = \cpf_{p,q}^{(2)}(\bsy{U})$, where $\bsy{U} = \{(u_i, v_j): 0 \leq i \leq p, 0 \leq j \leq q\}$ with $u_i=a$ and $n_j=b$ for all $i, j$.

    \item[(v)] $G$ is described by Theorem~\ref{maintheorem}(v)-(vi) and $\cpf(G) = \cpf_{p,q}^{(2)}(\bsy{U})$, where $\bsy{U} = \{(u_i, v_j): 0 \leq i \leq p, 0 \leq j \leq q\}$ with $\bsy{u} = (u_0, \ldots, u_{p-1})$ and $\bsy{v} = (v_0, \ldots, v_{q-1})$ given by Theorem \ref{thm:gaydarov and hopkins} according to the structures of $G_A$ and $G_B$, respectively.  (The values of $u_p$ and $v_q$ are not used as edge-weight of $D_{p,q}$ and are irrelevant. ) 
\end{enumerate}
\end{theorem}

\begin{remark}
Recall that in \Cref{sec:affine} we present graphs $G$ for which $\cpf_{p,q}^{(2)}(\bsy{U}) = \cpf(G)$ when $\bsy{U}$ is affine.
The affine $\bsy{U}$ with $c \geq 1$ 
described in \Cref{affinecasetheorem}
corresponds to  \Cref{maintheorem2}(iii), 
while the case with $c=0$ is included in 
\Cref{maintheorem2}(iv)-(v). 
\end{remark}

\Cref{maintheorem2} has the assumption $N_G(0) \cap A \neq \emptyset$. Exchanging $A$ and $B$ in $G$ corresponds to exchanging 
$u_{i,j}$ and $v_{i,j}$ in $\bsy{U}$. Hence, joining the cases in which $A$ and $B$ are switched in \Cref{maintheorem},  we obtain all possible $G$ such that $\cpf(G)=\cpf_{p,q}^{(2)}(\bsy{U})$ for some $\bsy{U}$. 
It remains to determine all the possible weight sets $\bsy{U}$. To this end, we notice that for a given graph $G$, there might be different weight sets $\bsy{U}$ that satisfy 
$\cpf(G)=\cpf_{p.q}^{(2)}(\bsy{U})$. See, for example, in Case 2 of Section~\ref{sec:affine},  we have $\cpf_{p.q}^{(2)}(\bsy{U})= \cpf_{p.q}^{(2)}(\bsy{U}')$.   


Recall that $\cpf_{p,q}^{(2)}(\bsy{U})=\cpf_{p,q}^{(2)}(\bsy{U}') $ if and only if $\mathcal{MI}_{p,q}^{(2)}(\bsy{U})=\mathcal{MI}_{p,q}^{(2)}(\bsy{U}')$.
In \Cref{sec:affine} we showed that  $\mathcal{MI}_{p,q}^{(2)}(\bsy{U})$ contains exactly the maximal elements in the set 
$\{(\bsy{a}(P), \bsy{b}(P)): \ P \text{ from $(0,0)$ to $(p,q)$}\}$,  where $(\bsy{a}(P), \bsy{b}(P))$ is defined in \eqref{eq:ab_P}. 
We will determine all the possible $\bsy{U}$  by analyzing the set $\mathcal{MI}_{p,q}^{(2)}(\bsy{U})$ on a case-by-case basis. 
The results are summarized in the following table. 

\begin{center}
\begin{tabular}{|c|c|c|} 
\hline 
  Case & size of $\mathcal{MI}_{p,q}^{(2)}(\bsy{U})$ &  Characterizing theorem  \\ \hline 
\Cref{maintheorem}(i)   &  2 &    \Cref{thm:U-case-1-2}(i)  \\ 
\Cref{maintheorem}(ii)  &  2 &    \Cref{thm:U-case-1-2}(ii)  \\ \Cref{maintheorem}(iii)  &  $\binom{p+q}{p}$ &    \Cref{thm:U-case3} \\ 
\Cref{maintheorem}(iv)-(vi)  &  1 &    \Cref{thm:U-case-4-6}  \\ \hline 
\end{tabular}
\end{center} 

Our first result is for Case (iii) of Theorem~\ref{maintheorem}.  

\begin{theorem} \label{thm:U-case3}
  Let $G$ be described by Theorem~\ref{maintheorem}(iii) and $\cpf(G) = \cpf_{p,q}^{(2)}(\bsy{U})$. Then $\bsy{U}$ is uniquely determined,  
   as described by Equation \eqref{Eq:affine-matrix}  in Theorem~\ref{maintheorem2}(iii). 
\end{theorem} 
\begin{proof}
   As seen in the proof of \Cref{lemma:c=c'}, when $c>0$, 
    $(\bsy{a}(P), \bsy{b}(P))$ and $(\bsy{a}(Q), \bsy{b}(Q))$ are non-comparable under $\preceq$ for any two distinct lattice paths $P, Q$ from $(0,0)$ to $(p,q)$. This implies that
     $\mathcal{MI}_{p,q}^{(2)}(\bsy{U})$   has $\binom{p+q}{p}$ elements, one for each lattice path $P$.  In particular, for any fixed index $i <p$,   there are $q+1$ different possible values, namely, $\{u_{i,j}: 0 \leq j \leq q\}$,  for the weights on the set of horizonal edges
     $E_i=\{ \text{ edge from} (i,j) \text{ to } (i+1, j):  0 \leq j \leq q\}$ in $D_{p,q}(\bsy{U})$. 
     
     Assume 
    there is another $\mathcal{U}'=\{ (u'_{i,j}, v'_{i,j}): 0\leq i \leq p, 0 \leq j \leq q\}$ such that $\mathcal{MI}_{p,q}^{(2)}(\bsy{U}')=\mathcal{MI}_{p,q}^{(2)}(\bsy{U})$. 
     Then in the digraph $D_{p,q}$ with weight  $\bsy{U}'$, weights on the edges in the set  $E_i$  must also have $q+1$ different values.  
    It forces $u_{i,j}=u'_{i,j}$ for all $j$. Similar argument shows $v_{i,j}=v'_{i,j}$, and hence, 
    $\bsy{U}'=\bsy{U}$.    
\end{proof}

Next we describe the set $\mathcal{MI}_{p,q}^{(2)}(\bsy{U})$ for the other cases of Theorem~\ref{maintheorem}.  
The proof is  a straightforward computation and hence omitted.

\begin{proposition} \label{prop:MI} 
  Assume the same conditions and the cases as in 
    Theorem~\ref{maintheorem}.  Let $\bsy{1}$ be the all $1$ vector of length $p+q$. 
     \begin{enumerate}
         \item[(i)]  $G$ is described by Theorem~\ref{maintheorem}(i) and $\cpf(G) = \cpf_{p,q}^{(2)}(\bsy{U})$ if and only if $\mathcal{MI}_{p,q}^{(2)}(\bsy{U})$ is the following set: 
         \begin{enumerate}
             \item[(a)]  If $p=1$, 
         $\mathcal{MI}_{p,q}^{(2)}(\bsy{U}) = \{ (a+b; b \dots, b)-\bsy{1}, (a; b, \dots, b, 2b)-\bsy{1}\}$. 
         
         If $p=2$,  $\mathcal{MI}_{p,q}^{(2)}(\bsy{U}) = \{ (a, a+b; b \dots, b)-\bsy{1}, (a, a; b, \dots, b, 2b)-\bsy{1}\}$. 

         \item[(b)] For $p > 2$, 
         $\mathcal{MI}_{p,q}^{(2)}(\bsy{U}) = \{ (a, \dots, a; a, \dots, a, 2a )-\bsy{1} , (a, \dots, a, 2a; a, \dots, a)-\bsy{1} \}. 
         $ 
\end{enumerate}

          \item[(ii)] $G$ is described by Theorem~\ref{maintheorem}(ii) and $\cpf(G) = \cpf_{p,q}^{(2)}(\bsy{U})$ if and only if          $$\mathcal{MI}_{p,q}^{(2)}(\bsy{U})=\{(a, a+b+c; c, \dots, c)-\bsy{1}, (a, a+b; c, \dots, c, 2c) -\bsy{1}\}.$$

          \item[(iii)] $G$ is described by Theorem~\ref{maintheorem}(iv)-(vi) and $\cpf(G) = \cpf_{p,q}^{(2)}(\bsy{U})$ if and only if  $\mathcal{MI}_{p,q}^{(2)}(\bsy{U})$  contains exactly one element $(\bsy{u}, \bsy{v})-\bsy{1}$, where 
          $\bsy{u} \in \{(a, \dots, a), (a, \dots, a, 2a), (a, a+b, \dots, a+(p-1)b)\}$ and $\bsy{v} \in \{(c, \dots, c), (c, \dots, c, 2c), (c, c+d, \dots, c+(q-1)d) \}. $ 
                  
     \end{enumerate}
\end{proposition}
\medskip

The following theorem constructs all $\bsy{U}$ when $\mathcal{MI}_{p,q}^{(2)}(\bsy{U})$ contains only one element, which solves the case that $G$ is described by Theorem~\ref{maintheorem}(iv)-(vi). 
\begin{theorem} \label{thm:U-case-4-6}
   Assume $\mathcal{MI}_{p,q}^{(2)}(\bsy{U})$ contains only one element $(\bsy{a},\bsy{b})$. Then the weight set $\bsy{U}$ can be constructed in the following steps. 
\begin{enumerate}
    \item Take a lattice path $P$ and put elements in $\bsy{a}$ as weights of horizonal edges of $P$, and elements in $\bsy{b}$ as weights  of vertical edges of $P$. 
    \item If the edge $\alpha=(i, j) \to (i+1, j)$ is in $P$, then any edge from $(i', j)$ to $(i'+1, j)$ has the same weight as $\alpha$, for $i'>i$. 
    \item If the edge $\beta=(i,j)\to (i, j+1)$ is in $P$, then any edge from $(i,j')$ to  $(i, j'+1)$ has the same weight as $\beta$, for $j' >j$. 
    \item The remaining edges can have any weight, as long as $u_{i,j} \leq u_{i',j'} $ and $v_{i,j} \leq v_{i',j'} $ when 
    $(i,j) \preceq (i',j')$. 
\end{enumerate}
\end{theorem}

For each of the remaining cases of \Cref{maintheorem},  there are two elements in
$\mathcal{MI}_{p,q}^{(2)}(\bsy{U})$.  The set of possible $\bsy{U}$ are described
in the next theorem. 

 \begin{theorem} \label{thm:U-case-1-2}
Consider the cases (i) and (ii) in 
 Theorem~\ref{maintheorem}.  
\begin{enumerate}
    \item[(i)] $G$ is described by Theorem~\ref{maintheorem}(i) and $\cpf(G) = \cpf_{p,q}^{(2)}(\bsy{U})$, where $\bsy{U}$
    is given by:  

    \begin{enumerate}
        \item If $p=1$, then $u_{0,q}=a+b, u_{0, q-1}=a$, $v_{1,q-1}=2b$, and  $v_{i,j}=b$ for $(i,j)\neq (1, q-1)$. Otherwise 
        $u_{0,0} \leq \cdots \leq u_{0, q-2}$ can have any values in $\{1, 2, \dots, a\}$. 
       \item If $p=2$,then $u_{1,q}=a+b$, $u_{0,q}=u_{0, q-1}=u_{1, q-1}=a$, $v_{p, q-1}=2b$, and $v_{i,j}=b$ for $(i,j)\neq (2, q-1)$. Other $u_{i,j}$'s can have any value
       in $\{1, 2, \dots, a\}$ as long as $u_{i,j} \leq u_{i',j'}$ whenever $(i,j) \preceq (i',j')$. 
       \item If $p \geq 2$,  then $u_{p-1,q} = v_{p,q-1} = 2a$, 
        $u_{p-1, q-1}=v_{p-1, q-1}=a$, and $u_{p,j} = v_{i,q} = a$ for all $i, j$. Then the remaining weights $\bsy{U'}=\{(u_{i,j}, v_{i,j}): 
        0 \leq i \leq p-1, 0 \leq j \leq q-1\}$ is obtained by taking any weights with $\mathcal{MI}_{p-1,q-1}^{(2)}(\bsy{U}')=\{ (a, \dots, a; a, \dots, a) \}$, as described in Theorem~\ref{thm:U-case-4-6}.       
    \end{enumerate}
    \item[(ii)] $G$ is described by Theorem~\ref{maintheorem}(ii) and  $\cpf(G) = \cpf_{p,q}^{(2)}(\bsy{U})$, where $p=2$. Then $\bsy{U}$
    is given by $v_{2,q-1} = 2c$, $v_{1, q-1}=c$, and $v_{2,j} = c$ for $j \neq q-1$; 
     $u_{0,q} =u_{0,q-1}= a$, $u_{1,q-1} = a+b$, $u_{1,q}=a+b+c$.
     For other edge weights, there are two possibilities. 
     \begin{enumerate}
         \item All $v_{i,j}=c$. 
     For other $u_{i,j}$, 
     we have $u_{0,j} \leq a$, $u_{1,j}\leq a+b$ for all $j < q-1$, and  $u_{i,j} \leq u_{i',j'}$ whenever $(i,j) \preceq (i',j')$. 
     \item  $u_{0, j}=a$ and $v_{1,j}=c$ for all $j$.  For other weights,
     $u_{1,0} \leq u_{1,1} \leq \cdots  \leq u_{1, q-2}$ can have any value in $\{a, \dots, a+b\}$, and $v_{0,0} \leq v_{0,1} \leq \cdots \leq v_{0,q-1}$ can have any values in $\{1, \dots, c\}$. 
      \end{enumerate}  
\end{enumerate}
\end{theorem} 
\begin{proof}
     We  present a detailed proof  for item (ii);  the cases in item (i) follow  similarly. 

     For item (ii), the maximal value at each position implies that $u_{0,q}=a, u_{1,q}=a+b+c$, $v_{2, q-1}=2c$ and $v_{2,j}=c$ for $j \leq q-2$.  Since $a+b+c$ and $2c$ cannot appear in the same element in $\mathcal{MI}_{p,q}^{(2)}(\bsy{U})$, we have 
     $u_{0, q-1}=a, u_{1,q-1}=a+b$,  and $v_{1,q-1}=c$. 

     Consider the lattice paths that correspond to elements in
     $\mathcal{MI}_{p,q}^{(2)}(\bsy{U})$. 
     \begin{enumerate}
         \item[(a)] If there is a lattice path $P$ such that 
         $P$ contains the edge $(0,0) \to (1,0)$ and $(\bsy{a}(P), \bsy{b}(P)) \in \mathcal{MI}_{p,q}^{(2)}(\bsy{U})$, 
         then $v_{0,0}=c$. This forces $v_{i,j}=c$ for all $i,j$.
         Now the two lattice paths $E^q N^2$ and $E^{q-1}NEN$ give the two  elements in $\mathcal{MI}_{p,q}^{(2)}(\bsy{U})$. Hence, other weights can be arbitrary as long as they satisfy the condition $u_{i,j} \leq u_{i',j'}$ whenever $(i,j) \preceq (i',j')$.

         \item[(b)]  If there is no lattice path $P$ such that $P$ contains the edge
         $(0,0) \to (1,0)$ and $(\bsy{a}(P), \bsy{b}(P)) \in \mathcal{MI}_{p,q}^{(2)}(\bsy{U})$.  Then $(\bsy{a}(Q), \bsy{b}(Q)) \in \mathcal{MI}(\bsy{U})$ for 
         the path $Q=NE^qN$.  This means $u_{0,0}=a$ and $v_{1,j}=c$ for all $j$.  It follows that $u_{0,j}=a$ for all $j$, as $u_{0,0} \leq u_{0,j} \leq u_{0, q}$.  
         Other weights can be arbitrary as long as they satisfy the condition that $u_{i,j} \leq u_{i',j'}$ whenever $(i,j) \preceq (i',j')$. In other words, $u_{1,0} \leq u_{1,1} \leq \cdots  \leq u_{1, q-2}$ can have any value in $\{a, \dots, a+b\}$, and $v_{0,0} \leq v_{0,1} \leq \cdots \leq v_{0,q-1}$ can have any values in $\{1, \dots, c\}$.       
         \end{enumerate}    
\end{proof}

Theorems~\ref{thm:U-case3}, \ref{thm:U-case-4-6} and \ref{thm:U-case-1-2}
together 
give a complete description of all the weight set $\bsy{U}$ such that there exists a graph $G$ with $\cpf(G)=\cpf_{p,q}^{(2)}(\bsy{U})$.

\let\oldbibitem\bibitem
\renewcommand{\bibitem}{\setlength{\itemsep}{0pt}\oldbibitem}

\bibliographystyle{amsplain}

\renewcommand{\bibname}{{\normalsize REFERENCES}}

\bibliography{reference}

\end{document}